\documentclass[10pt]{amsart}
\usepackage{amscd}
\usepackage{amssymb}
\usepackage[all]{xy}
\date{9 December 2006}
\title{Rigid Complexes via DG Algebras} 
\author{Amnon Yekutieli and James J. Zhang }
\address{A. Yekutieli: Department of  Mathematics 
Ben Gurion University, 
Be'er Sheva 84105, 
Israel}
\email{amyekut@math.bgu.ac.il}
\address{J.J. Zhang: Department of Mathematics, Box 354350,
University of Washington, Seattle, Washington 98195, USA}
\email{zhang@math.washington.edu}
\thanks{{\em Mathematics Subject Classification} 2000.
Primary: 18E30; Secondary: 18G10, 16E45, 18G15.}
\keywords{commutative rings, DG algebras,
derived categories, rigid complexes.}
\thanks{This research was supported by the US-Israel Binational
Science Foundation. The second author was partially supported by 
the US National Science Foundation.}

\newtheorem{thm}[equation]{Theorem}
\newtheorem{cor}[equation]{Corollary}
\newtheorem{prop}[equation]{Proposition}
\newtheorem{lem}[equation]{Lemma}
\theoremstyle{definition}
\newtheorem{dfn}[equation]{Definition}
\newtheorem{rem}[equation]{Remark}
\newtheorem{exa}[equation]{Example}

\newtheorem{conv}[equation]{Convention}

\numberwithin{equation}{section}

\newcommand{\iso}{\xrightarrow{\simeq}}
\newcommand{\inj}{\hookrightarrow}

\newcommand{\xar}{\xrightarrow}

\newcommand{\opn}{\operatorname}
\newcommand{\cat}[1]{\operatorname{\mathsf{#1}}}

\newcommand{\rmitem}[1]{\item[\text{\textup{(#1)}}]}
\newcommand{\mfrak}[1]{\mathfrak{#1}}
\newcommand{\mcal}[1]{\mathcal{#1}}
\newcommand{\msf}[1]{\mathsf{#1}}
\newcommand{\mbf}[1]{\mathbf{#1}}
\newcommand{\mrm}[1]{\mathrm{#1}}
\newcommand{\mbb}[1]{\mathbb{#1}}
\newcommand{\gfrac}[2]{\left[ \genfrac{}{}{0pt}{}{#1}{#2} \right]}

\newcommand{\tup}[1]{\textup{#1}}
\newcommand{\bsym}[1]{\boldsymbol{#1}}
\newcommand{\boplus}{\bigoplus\nolimits}

\newcommand{\til}[1]{\tilde{#1}}

\newcommand{\bra}[1]{\langle #1 \rangle}

\newcommand{\K}{\mbb{K} \hspace{0.05em}}
\renewcommand{\d}{\mathrm{d}}
\newcommand{\bwedge}{\bigwedge\nolimits}

\begin{document}

\begin{abstract}
Let $A$ be a commutative ring, $B$ a commutative $A$-algebra and 
$M$ a complex of $B$-modules. We begin by constructing the 
square $\opn{Sq}_{B / A} M$, which is also a complex of 
$B$-modules. The squaring operation is a quadratic functor, 
and its construction requires differential graded (DG) algebras.
If there exists an isomorphism $\rho : M \iso \opn{Sq}_{B / A} M$
then the pair $(M, \rho)$ is called a {\em rigid complex} over $B$ 
relative to $A$ (there are some finiteness conditions). 
There is an obvious notion of rigid morphism between rigid 
complexes. 

We establish several properties of rigid complexes, including 
their uniqueness, existence (under some extra hypothesis),
and formation of pullbacks 
$f^{\flat}(M, \rho)$ (resp.\ $f^{\sharp}(M, \rho)$)
along a finite (resp.\ essentially smooth) ring homomorphism
$f^* : B \to C$. 

In the subsequent paper \cite{YZ4} we consider {\em rigid dualizing 
complexes} over commutative rings, 
building on the results of the present paper. The project 
culminates in \cite{Ye4}, where we give a comprehensive version of 
Grothendieck duality for schemes. 

The idea of rigid complexes originates in noncommutative 
algebraic geometry, and is due to Van den Bergh \cite{VdB}.
\end{abstract}

\maketitle

%% ** section 0 **
\setcounter{section}{-1}
\section{Introduction}

{\rm Rigid dualizing complexes were invented by Van den Bergh 
\cite{VdB} in the context of noncommutative algebraic geometry. 
Since then the theory of rigid dualizing complexes was developed 
further by several people, with many applications in the areas of 
noncommutative algebra and noncommutative algebraic 
geometry. See the papers \cite{YZ1,YZ2,YZ3} and their references. 

The present paper is the first in a series of three papers in 
which we apply the rigidity technique to the areas of 
{\em commutative} algebra and algebraic geometry. Of course the 
commutative setup is ``contained in'', and is 
always much ``easier'' than the corresponding 
noncommutative setup. The point is that we want to get stronger 
and more general results when we restrict attention to commutative 
rings and schemes.

Let us remind the reader that in the papers cited above the 
noncommutative algebras in question were always assumed 
to be over some base 
{\em field} $\K$. The first generalization that comes to mind when 
considering a commutative theory of rigid dualizing complexes 
is to be able to work over {\em any commutative base ring $\K$}, 
or at least a sufficiently nice base ring, such as
$\K = \mbb{Z}$ for instance. This seemingly innocent 
generalization turns out to be quite hard. The present paper 
is devoted to solving it, by developing the technique of rigid 
complexes over commutative rings. 
The subsequent papers \cite{YZ4, Ye4}
will build on the foundations laid here, and will produce a 
comprehensive theory of dualizing complexes (including
duality for proper morphisms)
parallel to that of \cite{RD}, for finite type schemes over a 
regular finite dimensional commutative base ring $\K$. 

It should be noted that the rigidity theory developed here is 
essential even when dealing with commutative algebras (or schemes) 
over a base field $\K$. This is because it enables the concept of 
relative rigid complex, with respect to a $\K$-algebra homomorphism 
$A \to B$. 

For the rest of the introduction let us make the blanket 
assumption that all rings occurring are commutative. As usual, 
by an algebra over a ring $A$ we mean a ring $B$ together with a 
ring homomorphism $A \to B$. We denote by $\msf{D}(\cat{Mod} B)$
the derived category of $B$-modules. If $B$ is noetherian then we 
also consider the full subcategory 
$\msf{D}^{\mrm{b}}_{\mrm{f}}(\cat{Mod} B)$
consisting of bounded complexes with finitely generated 
cohomologies.

The difficulty in generalizing rigidity to an arbitrary base
ring $A$ shows up already in the definition of rigidity (Van den 
Bergh's rigidity equation). Let us first recall the 
original definition from \cite{VdB}. Suppose $A$ is a field,
$B$ is a noetherian 
$A$-algebra and $M \in \msf{D}^{\mrm{b}}_{\mrm{f}}(\cat{Mod} B)$.
A {\em rigidifying isomorphism} for $M$ is an isomorphism
\begin{equation} \label{eqn0.1}
\rho : M \iso \mrm{RHom}_{B \otimes_A B}(B, M \otimes_A M)
\end{equation}
in $\msf{D}(\cat{Mod} B)$. The pair $(M, \rho)$ is called a {\em 
rigid complex} over $B$ relative to $A$. 

Now suppose $A$ is any ring, and $B$ is a noetherian $A$-algebra 
which is {\em not flat}. Then in (\ref{eqn0.1}) we must replace 
$M \otimes_A M$ with $M \otimes^{\mrm{L}}_{A} M$; 
but $M \otimes^{\mrm{L}}_{A} M$ usually 
is not a well-defined object in
the derived category $\msf{D}(\cat{Mod} B \otimes_A B)$...

Our solution is to use differential graded (DG) algebras.
Take any ring $A$ and any $A$-algebra $B$. 
Choose some semi-free super-commutative non-positive DG 
algebra resolution $\til{B} \to B$ relative to $A$.
(All necessary facts about DG algebras are either reviewed or 
proved in Section 1 of the paper.)
Let $M \in \msf{D}(\cat{Mod} B)$. 
We prove that the complex
\[ \opn{Sq}_{B / A} M := 
\mrm{RHom}_{\til{B} \otimes_A \til{B}}(B, M \otimes^{\mrm{L}}_A M)
\in \msf{D}(\cat{Mod} B) \]
is independent of the resolution $\til{B} \to B$. Moreover, 
\[ \opn{Sq}_{B / A} : \msf{D}(\cat{Mod} B) \to \msf{D}(\cat{Mod} B)
\]
is a (nonlinear) functor. See Theorems \ref{thm1.1} and 
\ref{thm1.2}.

We are now ready to define rigidity. Let $A$ be a ring, let $B$ be 
a noetherian $A$-algebra and let 
$M \in \msf{D}^{\mrm{b}}_{\mrm{f}}(\cat{Mod} B)$. Assume $M$ has 
finite flat dimension over $A$. (This latter condition is 
automatically satisfied if $A$ is a regular ring.) Then, just like 
(\ref{eqn0.1}), a rigidifying isomorphism for $M$ is an 
isomorphism 
\[ \rho : M \iso \opn{Sq}_{B / A} M \]
in $\msf{D}(\cat{Mod} B)$, and the pair $(M, \rho)$ is called a 
rigid complex over $B$ relative to $A$. 

Suppose we are given two rigid complexes $(M, \rho_M)$ and
$(N, \rho_N)$ over $B$ relative to $A$. A morphism 
$\phi : M \to N$ in $\msf{D}(\cat{Mod} B)$ is said to be
{\em rigid} relative to $A$ if the diagram
\[ \UseTips \xymatrix @C=5ex @R=5ex {
M 
\ar[r]^(0.35){\rho_M}
\ar[d]_{\phi}
& \opn{Sq}_{B / A} M
\ar[d]^{\opn{Sq}_{B / A}(\phi)} \\
N
\ar[r]^(0.35){\rho_N}
& \opn{Sq}_{B / A} N
} \]
is commutative. The category of such rigid complexes and rigid 
morphisms is denoted by 
$\msf{D}^{\mrm{b}}_{\mrm{f}}(\cat{Mod} B)_{\mrm{rig} / A}$.

The following result about rigid complexes explains their name.

\begin{thm} \label{thm0.3}
Let $A$ be a ring, let $B$ be a noetherian $A$-algebra, and let
$(M, \rho) \in 
\msf{D}^{\mrm{b}}_{\mrm{f}}(\cat{Mod} B)_{\mrm{rig} / A}$.
Assume the canonical homomorphism 
$B \to \opn{Hom}_{\msf{D}(\cat{Mod} B)}(M, M)$
is bijective. Then the only automorphism of $(M, \rho)$
in 
$\msf{D}^{\mrm{b}}_{\mrm{f}}(\cat{Mod} B)_{\mrm{rig} / A}$
is the identity $\bsym{1}_M$. 
\end{thm}

This result is repeated as Theorem \ref{thm4.1} in the body of the 
paper. The proof boils down to the the following fact: 
let $\phi : M \to M$ be an isomorphism in 
$\msf{D}(\cat{Mod} B)$, so $\phi = b \bsym{1}_M$ for some 
invertible element $b \in B$. If $\phi$ is rigid then a 
calculation shows that $b^2 = b$; and hence $b = 1$. 

We find it convenient to denote ring homomorphisms by $f^*$ etc. 
Thus a ring homomorphism $f^* : A \to B$ corresponds to a 
morphism of schemes $f : \opn{Spec} B \to \opn{Spec} A$.

Let $A$ be a noetherian ring. 
Recall that an $A$-algebra $B$ is called 
essentially finite type if it is a localization of some finitely 
generated $A$-algebra. We say that $B$ is {\em essentially smooth}
(resp.\ {\em essentially \'etale}) over $A$ if it is essentially 
finite type and formally smooth (resp.\ formally \'etale).
These concepts are studied in Section 3. 

Let $A$ be a noetherian ring and $f ^* : A \to B$ an 
essentially smooth homomorphism. Then $\Omega^1_{B / A}$ is a 
finitely generated projective $B$-module. Let
$\opn{Spec} B = \coprod_i \opn{Spec} B_i$
be the decomposition into connected components, and for every $i$ 
let $n_i$ be the rank of $\Omega^1_{B_i / A}$. We define a functor
\[ f^{\sharp} : \msf{D}(\cat{Mod} A) \to \msf{D}(\cat{Mod} B) \]
by
\[ f^{\sharp} M := \bigoplus_i \, 
\Omega^{n_i}_{B_i / A}[n_i] \otimes_A M . \]

Recall that a ring homomorphism $f^* : A \to B$ is called finite 
if $B$ is a finitely generated $A$-module. Given such a finite 
homomorphism we define a functor
\[ f^{\flat} : \msf{D}(\cat{Mod} A) \to \msf{D}(\cat{Mod} B) \]
by
\[ f^{\flat} M := \opn{RHom}_A(B, M) . \]

\begin{thm} \label{thm0.1}
Let $A$ be a noetherian ring, let $B, C$ be essentially finite 
type $A$-algebras, let $f^* : B \to C$ be an $A$-algebra
homomorphism, and let $(M, \rho) \in 
\msf{D}^{\mrm{b}}_{\mrm{f}}(\cat{Mod} B)_{\mrm{rig} / A}$.
\begin{enumerate}
\item If $f^*$ is finite and $f^{\flat} M$ has finite flat 
dimension over $A$, then $f^{\flat} M$ has an induced rigidifying 
isomorphism 
\[ f^{\flat}(\rho) : f^{\flat} M  \iso \opn{Sq}_{C / A}
f^{\flat} M . \]
The rigid complex 
$f^{\flat}(\rho)(M, \rho) := 
\bigl( f^{\flat} M, f^{\flat}(\rho) \bigr)$
depends functorially on $(M, \rho)$ and on $f^*$. 
\item If $f^*$ is essentially smooth
then $f^{\sharp} M$ has an induced rigidifying 
isomorphism 
\[ f^{\sharp}(\rho) : f^{\sharp} M  \iso \opn{Sq}_{C / A}
f^{\sharp} M . \]
The rigid complex 
$f^{\sharp}(\rho)(M, \rho) := 
\bigl( f^{\sharp} M, f^{\sharp}(\rho) \bigr)$
depends functorially on $(M, \rho)$ and on $f^*$. 
\end{enumerate}
\end{thm}

This result is included in Theorems \ref{thm2.2} and \ref{thm2.4} 
in the body of the paper.

The next result is about tensor products of rigid complexes.

\begin{thm}
Let $A$ be a noetherian ring, let $B, C$ be essentially finite 
type $A$-algebras, let $f^* : B \to C$ be an $A$-algebra
homomorphism, and let $(M, \rho_M) \in 
\msf{D}^{\mrm{b}}_{\mrm{f}}(\cat{Mod} B)_{\mrm{rig} / A}$
and
$(N, \rho_N) \in 
\msf{D}^{\mrm{b}}_{\mrm{f}}(\cat{Mod} C)_{\mrm{rig} / B}$.
Assume that the canonical homomorphism 
$B \to \opn{Hom}_{\msf{D}(\cat{Mod} B)}(M, M)$
is bijective. Then the complex $M \otimes^{\mrm{L}}_{B} N$
is in $\msf{D}^{\mrm{b}}_{\mrm{f}}(\cat{Mod} C)$, it has finite 
flat dimension over $A$, and it has an 
induced rigidifying isomorphism 
\[ \rho_M \otimes \rho_N : M \otimes^{\mrm{L}}_{B} N \iso
\opn{Sq}_{C / A} (M \otimes^{\mrm{L}}_{B} N) . \]
\end{thm}

This is part of Corollary \ref{cor4.2}.

So far we said nothing about existence of rigid complexes. We can 
do so under some extra hypothesis. 

\begin{thm} \label{thm0.2}
Let $\K$ be a regular noetherian ring of finite Krull 
dimension, and let $A$ be an essentially finite type $\K$-algebra.
Then $A$ has a rigid complex $(R, \rho)$ relative to $\K$, 
such that the canonical homomorphism
$A \to \opn{Hom}_{\msf{D}(\cat{Mod} A)}(R, R)$
is bijective. 
\end{thm}

This theorem is repeated as Theorem \ref{thm6.4} in the body of the 
paper. Note that $R$ is nonzero (if $A \neq 0$), and moreover, by 
Theorem \ref{thm0.3}, the only rigid automorphism of 
$(R, \rho)$ is the identity. 
The construction of the rigid complex $(R, \rho)$ is actually 
quite easy (given Theorem \ref{thm0.1}).

The next theorem, which combines several results from our paper
\cite{YZ4}, gives a precise description of all rigid complexes in 
this situation. Recall that a {\em dualizing complex}
\cite{RD} is a complex 
$R \in  \msf{D}^{\mrm{b}}_{\mrm{f}}(\cat{Mod} A)$
which has finite injective dimension, and the canonical morphism 
$A \to \opn{RHom}_A(R, R)$ is an isomorphism.

\begin{thm}[\cite{YZ4}] \label{thm0.9}
Suppose $\K$ is a regular noetherian ring of finite Krull 
dimension. Let $A$ be an essentially finite type $\K$-algebra, such 
that $\opn{Spec} A$ is connected and nonempty. Then, up to 
isomorphism, the category 
$\msf{D}^{\mrm{b}}_{\mrm{f}}(\cat{Mod} A)_{\mrm{rig} / \K}$
has exactly two objects: the zero complex, and the 
nonzero rigid complex $(R, \rho)$ from Theorem \tup{\ref{thm0.2}}.
Moreover, $R$ is a dualizing complex over $A$. 
\end{thm}

We find Theorem \ref{thm0.9} very surprising, and we do not 
understand its significance yet. Also we do not know whether a 
similar result holds in the noncommutative setup. 

To end the introduction let us mention a couple of other
papers that approach 
Grothendieck duality in novel ways: \cite{Ne} and \cite{DGI}.

\medskip \noindent
\textbf{Acknowledgments.}
The authors wish to thank Bernhard Keller for his generous help 
with differential graded algebras. 
We also wish to thank Luchezar Avramov
and Michel Van den Bergh 
for useful discussions and valuable suggestions.
We are grateful to the referee for reading the paper carefully and 
making several corrections. 

\tableofcontents

%\newpage
%% ** section 1 **
\section{Differential Graded Algebras}
\label{sec1}

This section contains some technical material about differential 
graded algebras and their derived module categories. 
There is some overlap here with the papers \cite{FIJ},
\cite{Ke} and \cite{Be}.  

A graded algebra $A = \boplus_{i \in \mbb{Z}} A^i$
is said to be {\em super-commutative} if 
$a \cdot b = (-1)^{ij} b \cdot a$ for all $a \in A^i$ 
and $b \in A^j$,
and if $a \cdot a = 0$ whenever $i$ is odd. 
(Some authors call such a graded algebra strictly commutative.) 
The graded algebra 
$A$ is said to be {\em non-positive} if $A^i = 0$ for all 
$i > 0$. 

\begin{conv}
Throughout the paper all graded algebras are assumed by default
to be non-positive, super-commutative, associative and
unital. Homomorphisms of graded algebras are unital and have 
degree $0$. 
\end{conv}

Suppose $A \to B$ is a homomorphism of graded algebras. Then we 
say that $B$ is a graded $A$-algebra.

By {\em differential graded algebra} (or DG algebra) 
we mean a graded algebra 
$A = \bigoplus_{i \leq 0} A^i$,
together with a derivation 
$\mrm{d} : A \to A$ of degree $1$ satisfying 
$\mrm{d} \circ \mrm{d} = 0$. Note that the graded Leibniz rule 
holds: 
\[ \mrm{d}(a b) = \mrm{d}(a) b + (-1)^i a \mrm{d}(b) \]
for $a \in A^i$ and $b \in A^j$. 
The cohomology
$\mrm{H} A = \bigoplus_{i \leq 0} \mrm{H}^i A$
is then a graded algebra.

All rings in this paper are considered as DG algebras 
concentrated in degree $0$ (with zero differential); 
and in particular rings are assumed to be commutative. 

A DG algebra homomorphism $u : A \to B$ is a 
homomorphism of graded algebras
that commutes with $\mrm{d}$. It is a 
quasi-isomorphism if $\mrm{H}(u) : \mrm{H} A \to \mrm{H} B$ 
is an isomorphism (of graded algebras). 

If $A$ is a ring, $B$ is a DG algebra, and we are given a DG 
algebra homomorphism $A \to B$, then the differential $\d_B$ is 
necessarily $A$-linear. In this case we say that $B$ is a DG 
$A$-algebra. If $B = B^0$, i.e.\ $B$ is a ring too, then we say 
that $B$ is an $A$-algebra.

A differential graded (DG) $A$-module is a graded
(left) $A$-module $M = \boplus_{i \in \mbb{Z}} M^i$,
endowed with a degree $1$ $\mbb{Z}$-linear 
homomorphism $\mrm{d} : M \to M$ satisfying
$\mrm{d}(a m) = \mrm{d}(a) m + (-1)^i a \mrm{d}(m)$
for $a \in A^i$ and $m \in M^j$. Note that we can make $M$ into a 
right DG $A$-module by the rule 
$m a := (-1)^{i j} a m$ for $a \in A^i$ and $m \in M^j$. 
The category of DG $A$-modules is denoted by $\cat{DGMod} A$. 
It is an abelian category whose morphisms are degree $0$ 
$A$-linear homomorphisms commuting with the differentials.

There is a forgetful functor from DG algebras to graded algebras (it 
forgets the differential), and we denote it by
$A \mapsto \opn{und} A$. 
Likewise for $M \in \cat{DGMod} A$ we have
$\opn{und} M \in \cat{GrMod} (\opn{und} A)$, the category of graded 
$\opn{und} A$ -modules. If $A$ is a ring then a DG 
$A$-module is just a complex of $A$-modules. 

Given a graded algebra $A$ and two graded $A$-modules $M$ and $N$ 
let us write
\[ \opn{Hom}_{\mbb{Z}}(M, N)^i := \prod_{j \in \mbb{Z}}
\opn{Hom}_{\mbb{Z}}(M^j, N^{j+i}) , \]
the set of homogeneous $\mbb{Z}$-linear homomorphisms of degree $i$
from $M$ to $N$, and let
\[ \begin{aligned}
& \opn{Hom}_{A}(M, N)^i := \\
& \quad \{ \phi \in \opn{Hom}_{\mbb{Z}}(M, N)^i \mid
\phi(a m) = (-1)^{i j} a \phi(m) \text{ for all } 
a \in A^j \text{ and } m \in M \} .
\end{aligned} \]
Then 
\begin{equation} \label{eqn1.2}
\opn{Hom}_{A}(M, N) := \boplus_{i \in \mbb{Z}} 
\opn{Hom}_{A}(M, N)^i
\end{equation}
is a graded $A$-module, by the formula 
$(a \phi)(m) := a \phi(m) = (-1)^{i j} \phi(a m)$
for $a \in A^j$ and $\phi \in \opn{Hom}_{A}(M, N)^i$.
Cf. \cite[Chapter VI]{ML}. The set $\opn{Hom}_{A}(M, N)$ is related 
to the set of $A$-linear homomorphisms $M \to N$ 
as follows. Let's denote by $\opn{ungr}$ 
the functor forgetting the grading. Then the map
\[ \Phi : \opn{Hom}_{A}(M, N) \to 
\opn{Hom}_{\opn{ungr} A}(\opn{ungr} M, \opn{ungr}N) , \]
defined by
$\Phi(\phi)(m) := (-1)^{ij} \phi(m)$
for $\phi \in \opn{Hom}_{A}(M, N)^i$ and $m \in M^j$, is 
$\opn{ungr} A$ -linear, and $\Phi$ is bijective if $M$ 
is a finitely generated $A$-module.

For a DG algebra $A$ and two DG $A$-modules $M, N$ there is a 
differential $\mrm{d}$ on 
$\opn{Hom}_{\opn{und} A}(\opn{und} M, \opn{und} N)$, with formula
$\mrm{d}(\phi) := 
\mrm{d} \circ \phi - (-1)^{i} \phi \circ \mrm{d}$
for $\phi$ of degree $i$. 
The resulting DG $A$-module is denoted by 
$\opn{Hom}_{A}(M, N)$. 
Note that $\opn{Hom}_{\cat{DGMod} A}(M, N)$
coincides with the set of $0$-cocycles of 
$\opn{Hom}_{A}(M, N)$. 
Two homomorphisms 
$\phi_0, \phi_1 \in \opn{Hom}_{\cat{DGMod} A}(M, N)$
are  said to be {\em homotopic} if 
$\phi_0 - \phi_1 = \mrm{d}(\psi)$
for some $\psi \in \opn{Hom}_{A}(M, N)^{-1}$. The DG 
modules 
$M$ and $N$ are called {\em homotopy equivalent} if
there are homomorphisms 
$\phi : M \to N$ and $\psi : N \to M$ in $\cat{DGMod} A$ such that 
$\psi \circ \phi$ and $\phi \circ \psi$ are homotopic to the 
respective identity homomorphisms.

Suppose $A$ is a ring and $B, C$ are two DG $A$-algebras. Then 
$B \otimes_{A} C$ is also a DG $A$-algebra; the sign rule says 
that
\[ (b_1 \otimes c_1) \cdot (b_2 \otimes c_2) := 
(-1)^{i j} b_1 b_2 \otimes c_1 c_2 \]
for $c_1 \in C^j$ and $b_2 \in B^i$. 
The differential is of course
\[ \mrm{d}(b \otimes c) := \mrm{d}(b) \otimes c 
+ (-1)^i b \otimes \mrm{d}(c) \]
for $b \in B^i$. If $M \in \cat{DGMod} B$ and $N \in \cat{DGMod} C$
then $M \otimes_{A} N \in \cat{DGMod} B \otimes_{A} C$. 
If $N \in \cat{DGMod} B$ then $M \otimes_B N$, which is a quotient 
of $M \otimes_{A} N$, is a DG $B$-module. 

Let $A$ be a DG algebra.
Since $A$ is non-positive one has $\mrm{d}(A^0) = 0$; and 
therefore the differential $\mrm{d} : M^i \to M^{i+1}$ 
of any DG $A$-module $M$
is $A^0$-linear. This easily implies that the truncated objects
\begin{equation}
\begin{aligned}
\tau^{\geq i} M & := \big( \cdots 0 \to 
\opn{Coker}(M^{i-1} \to M^i) \to M^{i+1} \to \cdots \big) \\
& \text{and} \\
\tau^{\leq i} M & := \big( \cdots \to M^{i-1} \to
\opn{Ker}(M^{i} \to M^{i+1}) \to 0 \to \cdots \big) \\
\end{aligned}
\end{equation}
are DG $A$-modules. 

There is a derived category obtained from $\cat{DGMod} A$
by inverting the quasi-iso\-mor\-phisms, which we denote by
$\til{\msf{D}}(\cat{DGMod} A)$.
See \cite{Ke} for details. Note that
in case $A$ is a ring then 
$\cat{DGMod} A = \msf{C}(\cat{Mod} A)$,
the abelian category of complexes of $A$-modules, and 
$\til{\msf{D}}(\cat{DGMod} A) = \msf{D}(\cat{Mod} A)$,
the usual derived category of $A$-modules
(as in \cite{RD} or \cite{KS}).

In order to derive functors one has 
several useful devices. A DG $A$-module $P$ is called 
{\em K-projective} if for any acyclic DG $A$-module $N$ the 
DG module $\opn{Hom}_{A}(P, N)$ is acyclic.
(This name is due to Spaltenstein \cite{Sp}. Keller \cite{Ke} uses
the term ``property (P)'' to indicate K-projective DG modules,
and in \cite{AFH} the authors use ``homotopically projective''.
See also \cite{Hi}.)
Similarly one defines {\em K-injective} and {\em K-flat}
DG modules: $I$ is K-injective, and $F$ is K-flat, if 
$\opn{Hom}_{A}(N, I)$ and $F \otimes_A N$ are acyclic for 
all acyclic $N$. 
It is easy to see that any K-projective DG module is also K-flat.
Every two objects $M, N \in \cat{DGMod} A$
admit quasi-isomorphisms $P \to M$, $N \to I$ and $F \to M$,
with $P$ K-projective, $I$ K-injective and $F$ K-flat.
Then one defines
\[ \mrm{RHom}_{A}(M, N) := \mrm{Hom}_{A}(P, N) \cong
\mrm{Hom}_{A}(M, I) \in \til{\msf{D}}(\cat{DGMod} A) \]
and
\[ M \otimes^{\mrm{L}}_{A} N := F \otimes_{A} N \in
\til{\msf{D}}(\cat{DGMod} A) . \]

When $A$ is a usual algebra, any bounded above complex of projective 
(resp.\ flat) modules is K-projective (resp.\ K-flat). And any
bounded below complex of injective $A$-modules is K-injective. 
A single $A$-module $M$ is projective (resp.\ injective, 
resp.\ flat) iff it is K-projective (resp.\ K-injective, 
resp.\ K-flat) as DG $A$-module.

The following useful result is partly contained in 
\cite{Hi}, \cite{Ke} and \cite{KM}. 

\begin{prop} \label{prop1.1}
Let $A \to B$ be a quasi-isomorphism of DG algebras.
\begin{enumerate}
\item Given $M \in \til{\msf{D}}(\cat{DGMod} A)$ and 
$N \in \til{\msf{D}}(\cat{DGMod} B)$, the canonical morphisms 
$M \to B \otimes^{\mrm{L}}_{A} M$ and 
$B \otimes^{\mrm{L}}_{A} N \to N$ are both isomorphisms.
Hence the ``restriction of scalars'' functor
$\til{\msf{D}}(\cat{DGMod} B) \to \til{\msf{D}}(\cat{DGMod} A)$
is an equivalence.
\item Let $M, N \in \til{\msf{D}}(\cat{DGMod} B)$. Then there are
functorial isomorphisms 
$M \otimes^{\mrm{L}}_{B} N \cong M \otimes^{\mrm{L}}_{A} N$
and
$\mrm{RHom}_{B}(M, N) \cong \mrm{RHom}_{A}(M, N)$
in $\til{\msf{D}}(\cat{DGMod} A)$.
\end{enumerate}
\end{prop}

\begin{proof}
(1) Choose K-projective resolutions $P \to M$ and $Q \to N$
over $A$. Then $M \to B \otimes^{\mrm{L}}_{A} M$ 
becomes $P \cong A \otimes_A P \to B \otimes_A P$,
which is evidently a quasi-isomorphism. On the other hand
$B \otimes^{\mrm{L}}_{A} N \to N$ becomes
$B \otimes_A Q \to Q$; which is a quasi-isomorphism because so is 
$A \otimes_A Q \to B \otimes_A Q$.

\medskip \noindent
(2) Choose K-projective resolutions $P \to M$ and $Q \to N$
over $A$. We note that $B \otimes_A P$ and $B \otimes_A Q$ are
K-projective over $B$, and 
$B \otimes_A P \to M$, $B \otimes_A Q \to N$ are 
quasi-isomorphisms. Therefore we get isomorphisms in 
$\til{\msf{D}}(\cat{DGMod} A)$:
\[ M \otimes^{\mrm{L}}_{B} N = 
(B \otimes_A P) \otimes_B (B \otimes_A Q) \cong 
(B \otimes_A P) \otimes_A Q \cong P \otimes_A Q = 
M \otimes^{\mrm{L}}_{A} N . \]
The same resolutions give
\[ \mrm{RHom}_{B}(M, N) = \mrm{Hom}_{B}(B \otimes_A P, N) \cong
\mrm{Hom}_{A}(P, N) = \mrm{RHom}_{A}(M, N) . \]
\end{proof}

There is a structural characterization of K-projective DG
modules, which we shall review (since we shall elaborate on 
it later). This characterization works in steps. First one defines 
{\em semi-free} DG $A$-modules. A DG $A$-module $Q$ is called 
semi-free if there is a subset
$X \subset Q$ consisting of (nonzero) homogeneous elements, and an 
exhaustive non-negative increasing filtration 
$\{ F_{i} X \}_{i \in \mbb{Z}}$
of $X$ by subsets (i.e.\ $F_{-1} X = \emptyset$ and
$X = \bigcup F_i X$), such that
$\opn{und} Q$ is a free graded $\opn{und} A$ -module
with basis $X$, and for every $i$ one has
$\mrm{d}(F_i X) \subset \sum_{x \in F_{i - 1} X} A x$.
The set $X$ is called a {\em semi-basis} of $Q$. 
Note that $X$ is partitioned into
$X = \coprod_{i \in \mbb{Z}} X_i$,
where $X_i := X \cap Q^i$. 
We call such a set a {\em graded set}. 
Now a DG $A$-module $P$ is K-projective iff it is homotopy 
equivalent to a direct summand (in $\cat{DGMod} A$) of some 
semi-free DG module $Q$. 
See \cite{AFH} or \cite{Ke} for more details and for proofs. 

A {\em free \tup{(}super-commutative, non-positive\tup{)} 
graded $\mbb{Z}$-algebra} is a graded algebra of the following form. 
One starts with a graded set of variables 
$X = \coprod_{i \leq 0} X_{i}$; 
the elements of $X_i$ are the variables of 
degree $i$. Let 
$X_{\mrm{ev}} := \coprod_{i\ \text{even}} X_i$
and
$X_{\mrm{odd}} := \coprod_{i\ \text{odd}} X_i$.
Consider the free associative $\mbb{Z}$-algebra $\mbb{Z} \bra{X}$
on this set of variables. 
Let $I$ be the two-sided ideal of $\mbb{Z} \bra{X}$
generated by all elements of the form
$x y - (-1)^{i j} y x$ or
$z^2$, where
$x \in X_i$, $y \in X_j$, $z \in X_k$,
and $k$ is odd. The free super-commutative
graded $\mbb{Z}$-algebra on $X$ is the quotient
$\mbb{Z}[X] := \mbb{Z} \bra{X} / I$.
It is useful to note that
\[ \mbb{Z}[X] \cong \mbb{Z}[X_{\mrm{ev}}] 
\otimes_{\mbb{Z}} \mbb{Z}[X_{\mrm{odd}}] , \]
and that $\mbb{Z}[X_{\mrm{ev}}]$
is a commutative polynomial algebra, whereas
$\mbb{Z}[X_{\mrm{odd}}]$ is an exterior algebra.

\begin{dfn} \label{dfn1.2}
Suppose $A \to B$ is a homomorphism of DG algebras. 
$B$ is called a {\em semi-free} (super-commutative, non-positive)
{\em DG algebra relative to A} 
if there is a graded set $X = \coprod_{i \leq 0} X_i$, and 
an isomorphism of graded $\opn{und} A$ -algebras
\[ (\opn{und} A) \otimes_{\mbb{Z}} \mbb{Z}[X] \cong \opn{und} B  .\]
\end{dfn}

Observe that the DG algebra $B$ in the definition above, 
when regarded as a DG $A$-module, is semi-free with 
semi-basis consisting of the monomials in elements of $X$. Hence 
$B$ is also K-projective and K-flat as DG $A$-module. 

\begin{dfn}
Suppose $A$ and $B$ are DG algebras and $u: A \to B$ 
is a homomorphism of DG algebras. 
A {\em semi-free \tup{(}resp.\ K-projective, resp.\ K-flat\tup{)} 
DG algebra resolution of $B$ relative to A} is the data 
$A \xar{\til{u}} \til{B} \xar{v} B$,
where $\til{B}$ is a  DG 
$\K$-algebra, $\til{u}$ and $v$ are DG algebra homomorphisms, 
and the following conditions are satisfied:
\begin{enumerate}
\rmitem{i} $v \circ \til{u} = u$.
\rmitem{ii} $v$ is a quasi-isomorphism.
\rmitem{iii} $\til{u}$ makes $\til{B}$
into a semi-free DG algebra relative to $A$
\tup{(}resp.\ a K-projective DG $A$-module, 
resp.\ a K-flat DG $A$-module\tup{)}.
\end{enumerate}
We also say that $A \xar{\til{u}} \til{B} \xar{v} B$
is a {\em  semi-free \tup{(}resp.\ K-projective, resp.\ 
K-flat\tup{)}  DG algebra resolution} of $A \xar{u} B$. 
\[ \UseTips \xymatrix @C=5ex @R=3.5ex {
& \til{B}
\ar@{-->}[dr]^{v}
\\
A 
\ar@{-->}[ur]^{\til{u}} 
\ar[rr]^{u} 
& & B 
} \]
\end{dfn}

\begin{prop} \label{prop1.2}
Let $A$ and $B$ be  DG algebras, 
and let $u: A \to B$ be a DG algebra homomorphism. 
\begin{enumerate}
\item There exists a semi-free DG algebra resolution
$A \xar{\til{u}} \til{B} \xar{v} B$
of $A \xar{u} B$.
\item Moreover, if $\mrm{H} A$ is a noetherian algebra and 
$\mrm{H} B$ is a finitely generated $\mrm{H} A$ -algebra, 
then we can choose the semi-free DG algebra $\til{B}$ 
in part \tup{(1)} such that 
$\opn{und} \til{B} \cong (\opn{und} A) \otimes_{\mbb{Z}} \mbb{Z}[X]$
as graded $\opn{und} A$ -algebras, 
where the graded set $X = \coprod_{i \leq 0} X_i$
has finite graded components $X_{i}$. 
\item If $\mrm{H} A$ is a noetherian algebra, $B$ is a ring, 
and $B = \mrm{H}^0 B$ is a finitely generated 
$\mrm{H}^0 A$ -module, 
then there exists a K-projective DG algebra resolution
$A \to \til{B} \to B$ of $A \to B$,
such that 
$\opn{und} \til{B} \cong \boplus_{i = -\infty}^{0} 
\opn{und} A [-i]^{\mu_i}$ as graded $\opn{und} A$ -modules,
and the multiplicities $\mu_i$ are finite. 
\end{enumerate}
\end{prop}

\begin{proof}
(1) We shall construct $\til{B}$ as the union of an increasing 
sequence of DG algebras
$F_0 \til{B} \subset F_1 \til{B} \subset \cdots$, 
which will be defined recursively. 
At the same time we shall construct an increasing sequence 
of DG algebra homomorphisms
$A \to F_i \til{B} \xar{v_i} B$,
and an increasing sequence of graded sets 
$F_i X \subset F_i \til{B}$. The homomorphism $v$ will be the 
union of the $v_i$, and the graded set $X = \coprod_{j \leq 0} X_j$
will be the union of the sets $F_i X$.
For every $i$ the following conditions will hold:
\begin{enumerate}
\rmitem{i} $\mrm{H}(v_i) : \mrm{H}(F_i \til{B}) \to \mrm{H} B$
is surjective in degrees $\geq -i$.
\rmitem{ii} $\mrm{H}(v_i) : \mrm{H}(F_i \til{B}) \to \mrm{H} B$
is bijective in degrees $\geq -i + 1$.
\rmitem{iii} $F_i \til{B} = A[F_i X]$, 
$\mrm{d}(F_i X) \subset F_{i-1} \til{B}$ and
$\opn{und} F_i \til{B} \cong (\opn{und} A) \otimes_{\mbb{Z}}
\mbb{Z}[F_i X]$.
\end{enumerate}

We start by choosing a set of elements of $B^0$ that 
generate $\mrm{H}^0 B$ as $\mrm{H}^0 A$ -algebra. 
This gives us a set $X_0$ of elements of degree $0$ with a function 
$v_0 : X_0 \to B^0$. Consider the DG algebra
$\mbb{Z}[X_0]$ with zero differential; and define
$F_0 \til{B} := A \otimes_{\mbb{Z}} \mbb{Z}[X_0]$. 
Also define $F_0 X := X_0$. 
We get a DG algebra homomorphism
$v_0 : F_0 \til{B} \to B$, and conditions (i)-(iii) hold for 
$i = 0$.

Now assume $i \geq 0$, and that for every $j \leq i$
we have DG algebra homomorphisms
$v_j : F_j \til{B} \to B$ and graded sets $F_j X$ 
satisfying conditions (i)-(iii). 
We will construct $F_{i+1} \til{B}$ etc.

Choose a set $Y_{i + 1}'$ of elements (of degree $-i - 1$)
and a function $v_{i+1} : Y_{i+1}' \to B^{-i - 1}$ such that
$\{ v_{i+1}(y) \mid y \in Y_{i+1}' \}$
is a set of cocycles that generates
$\mrm{H}^{-i - 1} B$ as $\mrm{H}^0 A$-module. For 
$y \in Y_{i+1}'$ define $\mrm{d}(y) := 0$. 

Next let
\[ J_{i+1} := \{ b \in (F_i \til{B})^{-i} \mid
\mrm{d}(b) = 0 \text{ and } 
\mrm{H}^{-i}(v_i)(b) = 0 \} . \]
Choose a set $Y_{i + 1}''$ of elements (of degree $-i-1$)
and a function 
$\mrm{d} : Y_{i+1}'' \to J_{i+1}$ such that
$\{ \mrm{d}(y) \mid y \in Y_{i+1}'' \}$
is a set of elements whose images in $\mrm{H}^{-i} F_i \til{B}$
generate
$\opn{Ker} \bigl( \mrm{H}^{-i}(v_i) : 
\mrm{H}^{-i} F_i \til{B} \to \mrm{H}^{-i} B \bigr)$
as $\mrm{H}^0 A$-module. Let $y \in Y_{i+1}''$. By definition 
$v_i(\mrm{d}(y)) = \mrm{d}(b)$
for some $b \in B^{-i}$; and we define $v_{i+1}(y) := b$.  

Let $Y_{i+1} := Y_{i+1}' \sqcup Y_{i+1}''$
and $F_{i+1} X := F_i X \sqcup Y_{i+1}$. 
Define the DG algebra $F_{i+1} \til{B}$ to be
\[ F_{i+1} \til{B} := F_{i} \til{B} 
\otimes_{\mbb{Z}} \mbb{Z}[Y_{i+1}] \]
with differential $\mrm{d}$ extending the differential of 
$F_{i} \til{B}$ and the function
$\mrm{d} : Y_{i+1} \to F_{i} \til{B}$
defined above.

\medskip \noindent
(2) This is because at each step in (1) the sets $Y_i$ can be chosen 
to be finite.

\medskip \noindent
(3) Choose elements $b_1, \ldots, b_m \in B$ that generate it as 
$A^0$-algebra. Since each $b_i$ is integral over $A^0$, there is 
some monic polynomial 
$p_i(y) \in A^0[y]$ such that $p_i(b_i) = 0$. Let 
$y_1, \ldots, y_m$ be distinct variables of degree $0$. Define
$Y_0 := \{ y_1, \ldots, y_m \}$ and
$B^{\dag} := A^0[Y_0] / \big( p_1(y_1), \ldots, p_m(y_m) \big)$.
This is an $A^0$-algebra, which is a free module 
of finite rank. Let $v_0 : B^{\dag} \to B$ be the 
surjective $A^0$-algebra homomorphism $y_i \mapsto b_i$. Define
$F_0 \til{B} := A \otimes_{A^0} B^{\dag}$ and
$F_0 X := \emptyset$. Then conditions (i)-(ii) hold for $i = 0$, 
as well as condition (iii') below.
\begin{enumerate}
\rmitem{iii'} $F_i \til{B} = A[Y_0 \cup F_i X]$, 
$\mrm{d}(F_i X) \subset F_{i-1} \til{B}$ and
\[ \opn{und} F_i \til{B} \cong (\opn{und} A) \otimes_{A^0}
A^0[F_i X] \otimes_{A^0} B^{\dag} . \] 
\end{enumerate}

For $i \geq 1$ the proof proceeds as in part (1), but always using 
condition (iii') instead of (iii).
\end{proof}

\begin{prop} \label{prop1.3}
Suppose we are given three  
DG algebras  $\til{A}, \til{B}, \til{B}'$;
a ring $B$; and five DG algebra homomorphisms
$u, \til{u}, \til{u}', v, v'$ such that the first
diagram below is commutative. Assume that $v'$ is a 
quasi-isomorphism, and $\til{B}$ is semi-free DG algebra 
relative to $\til{A}$. Then there exists a DG algebra 
homomorphism $w : \til{B} \to \til{B}'$
such that the second diagram below is commutative. 
\[ \begin{array}{cc}
\UseTips \xymatrix @C=5ex @R=3.5ex {
& \til{B}
\ar[dr]^{v}
\\
\til{A}
\ar[ur]^{\til{u}} 
\ar[dr]_{\til{u}'} 
\ar[rr]^{u} 
& & B \\
& \til{B}' 
\ar[ur]_{v'} }
\qquad \quad 
& 
%
% second diagram
%
\UseTips \xymatrix @C=5ex @R3.5ex {
& \til{B}
\ar[dr]^{v}
\ar@{-->}[dd]^{w}
\\
\til{A}
\ar[ur]^{\til{u}} 
\ar[dr]_{\til{u}'} 
& & B \\
& \til{B}'
\ar[ur]_{v'} } 
\end{array}
\] 
\end{prop}

\begin{proof}
By definition there is a graded set 
$X = \coprod_{i \leq 0} X_{i}$ such that 
$\opn{und} \til{B} \cong$ \linebreak
$(\opn{und} \til{A}) \otimes_{\mbb{Z}} \mbb{Z}[X]$.
Let's define $F_i X := \bigcup_{j \geq -i} X_j$ and
$F_i \til{B} := \til{A}[F_i X] \subset \til{B}$. 
We shall define a compatible sequence of DG algebra homomorphisms 
$w_i : F_i \til{B} \to \til{B}'$,
whose union will be called $w$. 

For $i = 0$ we note that 
$v' : \til{B}'^0 \to B$ is surjective. Hence there is a function
$w_0 : X_0 \to \til{B}'^0$
such that $v'(w_0(x)) = v(x)$ for every $x \in X_0$. Since
$F_0 \til{B} \cong \til{A} \otimes_{\mbb{Z}} \mbb{Z}[X_{0}]$ 
and 
$\mrm{d}(w_0(X_0)) = 0$ we can extend the function 
$w_0$ uniquely to a DG algebra homomorphism
$w_0 : F_0 \til{B} \to \til{B}$
such that $w_0 \circ \til{u} = \til{u}'$. 

Now assume that $i \geq 0$ and 
$w_i : F_i \til{B} \to \til{B}'$
has been defined. Let $Y_{i+1} := F_{i+1} X - F_i X$. This is a 
set of degree $-i-1$ elements. Take any $y \in Y_{i+1}$. Then
$\mrm{d}(y) \in (F_i \til{B})^{-i}$, and we let
$b := w_i(\mrm{d}(y)) \in \til{B}'^{-i}$.
Because 
$\mrm{H} \til{B} \cong \mrm{H} \til{B}' = B$
there exists an element
$c \in  \til{B}'^{-i-1}$
such that $\mrm{d}(c) = b$. We now define
$w_{i+1}(y) := c$. 
The function $w_{i+1} : Y_{i+1} \to \til{B}'^{-i-1}$
extends to a unique DG algebra homomorphism
$w_{i+1} : F_{i+1} \til{B} \to \til{B}'$
such that $w_{i+1}|_{F_{i} \til{B}} = w_{i}$. 
\end{proof}

A homomorphism $A \to A'$ between two rings
is called a {\em localization} if it induces an isomorphism
$S^{-1} A \iso A'$ for some multiplicatively closed subset
$S \subset A$. We then say that $A'$ is a localization of $A$.

Suppose $A$ is a noetherian ring. An $A$-algebra 
$B$ is called {\em essentially of finite type} if 
$B$ is a localization of some finitely generated $A$-algebra.
Such an algebra $B$ is noetherian. If $C$ is an essentially 
finite type $B$-algebra then it is an essentially finite type 
$A$-algebra.

\begin{prop} \label{prop1.4}
Let $A$ be a noetherian ring, and
let $B$ be an essentially finite type $A$-algebra. 
Then there is a DG algebra quasi-isomorphism $\til{B} \to B$ such 
that $\til{B}^0$ is an essentially finite type $A$-algebra, 
and each $\til{B}^i$ is a finitely 
generated $\til{B}^0$-module and a flat $A$-module. In particular 
$\til{B}$ is a K-flat DG $A$-module. 
\end{prop}

\begin{proof}
Pick a finitely generated $A$-algebra $B_{\mrm{f}}$ such that 
$S^{-1} B_{\mrm{f}} \cong B$ for some multiplicatively closed 
subset $S \subset B_{\mrm{f}}$. According to Proposition 
\ref{prop1.2}(2) we can find a semi-free DG algebra resolution
$\til{B}_{\mrm{f}} \to B_{\mrm{f}}$ where $\til{B}_{\mrm{f}}$ has 
finitely many algebra generators in each degree. Let
$\til{S} \subset \til{B}_{\mrm{f}}^0$ be the pre-image of $S$
under the surjection $\til{B}_{\mrm{f}}^0 \to B_{\mrm{f}}$. 
Now define
$\til{B} := (\til{S}^{-1} \til{B}_{\mrm{f}}^0)
\otimes_{\til{B}_{\mrm{f}}^0} \til{B}_{\mrm{f}}$. 
\end{proof}

\begin{cor} \label{cor1.1}
Let $B$ be an essentially finite type 
$A$-algebra, and let $\til{B} \to B$ be any K-flat DG algebra 
resolution relative to $A$. Then 
$\mrm{H}^0(\til{B} \otimes_{A} \til{B})$
is an essentially finite type $A$-algebra, and each 
$\mrm{H}^i(\til{B} \otimes_{A} \til{B})$
is a finitely generated 
$\mrm{H}^0(\til{B} \otimes_{A} \til{B})$-module.
\end{cor}

\begin{proof}
Using Proposition \ref{prop1.3}, and passing via a semi-free 
DG algebra resolution, we can replace the given 
resolution $\til{B} \to B$ by another one satisfying the 
finiteness conditions in Proposition \ref{prop1.4}. Now the 
assertion is clear. 
\end{proof}

Let $M$ be a graded module. The {\em amplitude} 
$\opn{amp} M \in \mbb{N} \cup \{ \infty \}$ 
is defined as follows: if $M$ is nonzero and bounded then
\[ \opn{amp} M := 
\max \{ i \in \mbb{Z} \mid M^i \neq 0 \} - 
\min \{ i \in \mbb{Z} \mid M^i \neq 0 \} . \]
If $M$ is unbounded then $\opn{amp} M := \infty$, and 
if $M = 0$ then $\opn{amp} M := 0 $.

Now let $A$ be a DG algebra with $\mrm{H} A$ bounded, and let $M$ 
be a DG $A$-module. 
The {\em flat dimension} 
$\opn{flat{.}dim}_A M \in \mbb{N} \cup \{ \infty \}$ 
is defined as follows: 
\[ \begin{aligned}
& \opn{flat{.}dim}_A M := \\
& \qquad \inf\, \bigl\{ d \in \mbb{N} \mid 
\opn{amp} \mrm{H} (M \otimes^{\mrm{L}}_{A} N) \leq d +
\opn{amp} \mrm{H} N \text{ for all } N \in \cat{DGMod} A \bigr\} .
\end{aligned} \] 
Observe that $M$ has finite flat dimension if and only if the 
functor $M \otimes^{\mrm{L}}_{A} -$ is way out on both sides, in 
the sense of \cite[Section I.7]{RD}.
Also, if $M$ has finite flat dimension then $\mrm{H} M$ is bounded
(take $N := A$). 
Similarly one can define the projective dimension 
$\opn{proj{.}dim}_A M$ of a DG 
$A$-module $M$, by considering the amplitude of 
$\mrm{H}\, \mrm{RHom}_{A}(M, N)$. 
For a ring $A$ and a usual
module $M$ the dimensions defined above
coincide with the usual ones.

\begin{rem}
Our definition of flat dimension and projective dimension for DG 
modules differs from that of other authors, e.g.\ \cite{FIJ}.
\end{rem}

\begin{prop} \label{prop2.4}
Let $A$ be a ring, let $B$ and $C$ be DG $A$-algebras,  
$L \in \cat{DGMod} B$, $M \in$ 
$\cat{DGMod} B \otimes_{A} C$ and
$N \in \cat{DGMod} C$. 
There exists a functorial morphism
\[ \psi : 
\mrm{RHom}_{B}(L, M) \otimes^{\mrm{L}}_{C} N \to
\mrm{RHom}_{B}(L, M \otimes^{\mrm{L}}_{C} N) \]
in $\til{\msf{D}}(\cat{DGMod} B \otimes_{A} C)$.
If conditions \tup{(i)}, \tup{(ii)}, and \tup{(iii)}
below hold, then the morphism $\psi$ is an isomorphism.
\begin{enumerate}
\rmitem{i} $\mrm{H}^0 B$ is noetherian, $\mrm{H} L$ is 
bounded above, and each of the $\mrm{H}^0 B$ -modules 
$\mrm{H}^i B$ and $\mrm{H}^i L$ are finitely generated.
\rmitem{ii} $\mrm{H} M$ is bounded below, and $\mrm{H} N$ 
is bounded. 
\rmitem{iii} Either \tup{(a)}, \tup{(b)} or \tup{(c)} is satisfied:
\begin{enumerate}
\rmitem{a} $\mrm{H}^i B = 0$ for all $i \neq 0$, and 
$L$ has finite projective dimension over $B$. 
\rmitem{b} $\mrm{H}^i C = 0$ for all $i \neq 0$, and 
$N$ has finite flat dimension over $C$. 
\rmitem{c} $\mrm{H}^i C = 0$ for all $i \neq 0$, 
$\mrm{H}^0 C$ is noetherian, each
$\mrm{H}^i N$ is a finitely generated module over $\mrm{H}^0 C$,
the canonical morphism 
$C \to \mrm{RHom}_{C}(N, N)$ is an isomorphism, 
and both $M$ and $\mrm{RHom}_{B}(L, M)$ have finite 
flat dimension over $C$.
 \end{enumerate}
\end{enumerate}
\end{prop}

\begin{proof}
The proof is in five steps. 

\medskip \noindent
Step 1. To define $\psi$ we may choose a K-projective resolution 
$P \to L$ over $B$, and a K-flat resolution $Q \to N$ over $C$. 
There an obvious homomorphism of DG $B \otimes_{A} C$
-modules
\[ \psi_{P, Q} : 
\mrm{Hom}_{B}(P, M) \otimes_{C} Q \to
\mrm{Hom}_{B}(P, M \otimes_{C} Q) . \]
In the derived category this represents $\psi$. 

\medskip \noindent Step 2.
To prove that $\psi$ is an isomorphism (or equivalently that
$\psi_{P, Q}$ is a quasi-iso\-morphism)
we may forget the 
$B \otimes_{A} C$ -module structures, and consider $\psi$ as a 
morphism in $\msf{D}(\cat{Mod} A)$. 
Now by Proposition \ref{prop1.1}(2) we can replace $B$ and $C$ by 
quasi-isomorphic DG $A$-algebras. Thus we may assume both $B$ and 
$C$ are semi-free as DG $A$-modules. 

\medskip \noindent Step 3.
Let's suppose that condition (iii.a) holds. So 
$B \to \mrm{H}^0 B$ is a quasi-isomorphism. Since $C$ is K-flat 
over $A$ it follows that 
$B \otimes_{A} C \to \mrm{H}^0 B \otimes_{A} C$
is also a quasi-isomorphism. By Proposition \ref{prop1.1}
we can assume that 
$L \in \cat{DGMod} \mrm{H}^0 B$ and
$M \in \cat{DGMod} (\mrm{H}^0 B \otimes_{A} C)$,
and that $L$ has finite projective dimension over 
$\mrm{H}^0 B$. So we may replace $B$ with $\mrm{H}^0 B$, and thus
assume that $B$ is a noetherian algebra.

Now choose a resolution 
$P \to L$, where $P$ is a bounded complex of finitely generated 
projective $B$-modules. 
Take any K-flat resolution $Q \to N$ over $C$. 
Then the homomorphism $\psi_{P, Q}$ is actually bijective.

\medskip \noindent Step 4.
Let's assume condition (iii.b) holds. As in step 3
we can suppose that $C = C^0$. 
Choose a bounded resolution $Q \to N$ by flat $C$-modules.
By replacing $M$ with the truncation 
$\tau^{\geq j_0} M$ for some $j_0 \ll 0$, 
we may assume $M$ is bounded below.
According to \cite[Theorem 9.2.7]{AFH} we can find a semi-free 
resolution $P \to L$ over $B$ such that
$\opn{und} P \cong \boplus_{i = -\infty}^{i_1}
\opn{und} B[-i]^{\mu_i}$
with all the multiplicities $\mu_i$ finite. 
Because the $\mu_i$ are finite, $M$ is bounded below 
and $Q$ is bounded, 
the homomorphism $\psi_{P, Q}$ is bijective.

\medskip \noindent Step 5.
Finally we consider condition (iii.c). We can assume that $C = C^0$ 
is noetherian. Since 
$N \in \msf{D}^{\mrm{b}}_{\mrm{f}}(\cat{Mod} C)$
and
$\mrm{RHom}_{C}(N, N) \cong C$ 
we see that the support of $N$ is $\opn{Spec} C$. By Lemma 
\ref{lem9.5} below we conclude that 
$N$ generates 
$\msf{D}^{+}(\cat{Mod} C)$. Let
\[ \psi' : \mrm{RHom}_{C} \big( N, 
\mrm{RHom}_{B}(L, M) \otimes^{\mrm{L}}_{C} N \big) \to
\mrm{RHom}_{C} \big( N,
\mrm{RHom}_{B}(L, M \otimes^{\mrm{L}}_{C} N) \big) \]
be the morphism obtained from $\psi$ by applying the functor
$\mrm{RHom}_{C}(N, -)$. 
Since $\psi$ is a morphism in
$\msf{D}^{+}(\cat{Mod} C)$, 
in order to prove it is an isomorphism 
it suffices to prove that $\psi'$ is an isomorphism. 

Consider the commutative diagram of morphisms in
$\msf{D}(\cat{Mod} C)$
\[ \UseTips \xymatrix @C=7ex @R=4ex {
\mrm{RHom}_{B}(L, M) 
\ar[r]^(0.35){\alpha}
\ar[dr]!L^(0.6){\beta}
& \mrm{RHom}_{C} \big( N,
\mrm{RHom}_{B}(L, M) \otimes^{\mrm{L}}_{C} N \big) 
\ar[d]^{\psi'} \\
& \mrm{RHom}_{C} \big( N,
\mrm{RHom}_{B}(L, M \otimes^{\mrm{L}}_{C} N) \big) 
} \]
in which $\alpha$ and $\beta$ are the obvious morphisms. We shall 
prove that $\alpha$ and $\beta$ are isomorphisms. 
The complex $\mrm{RHom}_{B}(L, M)$ has finite flat dimension over 
$C$, so using the proposition with condition (iii.b), which we 
already proved, we have
\[ \begin{aligned}
\mrm{RHom}_{C} \big( N, 
\mrm{RHom}_{B}(L, M) \otimes^{\mrm{L}}_{C} N \big) & \cong
\mrm{RHom}_{C}(N, N) \otimes^{\mrm{L}}_{C} 
\mrm{RHom}_{B}(L, M) \\
& \cong \mrm{RHom}_{B}(L, M) . 
\end{aligned} \]
But the composed isomorphism is precisely $\alpha$. 
On the other hand the complex $M$ has finite flat dimension over 
$C$, so using the proposition with condition (iii.b) once more
(for the isomorphism marked $\diamond$), we have
\[ \begin{aligned}
& \mrm{RHom}_{C} \big( N, 
\mrm{RHom}_{B}(L, M \otimes^{\mrm{L}}_{C} N) \big) \cong
\mrm{RHom}_{B} \big( L, 
\mrm{RHom}_{C}(N, M \otimes^{\mrm{L}}_{C} N) \big) \\
& \quad \cong^{\diamond}
\mrm{RHom}_{B} \big( L, M \otimes^{\mrm{L}}_{C}
\mrm{RHom}_{C}(N, N) \big) 
\cong \mrm{RHom}_{B}(L, M) . 
\end{aligned} \]
Here the composed isomorphism is $\beta$.
\end{proof}

Let $C$ be a noetherian ring. Recall that given a complex 
$N \in \msf{D}^{\mrm{b}}_{\mrm{f}}(\cat{Mod} C)$
its support is defined to be
$\bigcup_i \opn{Supp} \mrm{H}^i N \subset \opn{Spec} C$. 
The complex $N$ is said to generate 
$\msf{D}^{+}(\cat{Mod} C)$ if
for any nonzero object $M \in \msf{D}^{+}(\cat{Mod} C)$
one has
$\mrm{RHom}_{C}(N, M) \neq 0$.

\begin{lem} \label{lem9.5}
Suppose $C$ is a noetherian ring and  
$N \in \msf{D}^{\mrm{b}}_{\mrm{f}}(\cat{Mod} C)$ 
is a complex whose support is $\opn{Spec} C$. Then $N$ 
generates $\msf{D}^{+}(\cat{Mod} C)$.
\end{lem}

\begin{proof}
Suppose $M$ is a nonzero object in 
$\msf{D}^{+}(\cat{Mod} C)$. We have to prove that 
$\mrm{RHom}_{C}(N, M) \neq 0$.
Let 
$i_0 := \opn{min} \{ i \in \mbb{Z} \mid \mrm{H}^i M \neq 0 \}$,
and choose a nonzero finitely generated submodule
$M' \subset \mrm{H}^{i_0} M$.
Let $\mfrak{p}$ be a minimal prime ideal in the support of 
$M'$; so that $M'_{\mfrak{p}} := B_{\mfrak{p}} \otimes_C M'$ 
is a nonzero finite length module 
over the local ring $C_{\mfrak{p}}$. Now 
$N_{\mfrak{p}}$ is a nonzero object of
$\msf{D}^{\mrm{b}}_{\mrm{f}}(\cat{Mod} C_{\mfrak{p}})$.
Let 
$j_1 := \opn{max} \{ j \in \mbb{Z} \mid \mrm{H}^j 
N_{\mfrak{p}} \neq 0 \}$. 
Since $\mrm{H}^{j_1} N_{\mfrak{p}}$ is a nonzero finitely generated 
$B_{\mfrak{p}}$-module, there exists a nonzero homomorphism
$\phi : \mrm{H}^{j_1} N_{\mfrak{p}} \to M'_{\mfrak{p}}$.
This $\phi$ can be interpreted as a nonzero element of
$\opn{Ext}_{C_{\mfrak{p}}}^{i_0 - j_1}(N_{\mfrak{p}},
M_{\mfrak{p}})$. 

Finally, by Proposition \ref{prop2.4} under its condition (iii.b)
-- whose proof does not rely on this lemma --
we have an isomorphism
\[ \opn{Ext}_{C_{\mfrak{p}}}^{i_0 - j_1}(N_{\mfrak{p}},
M_{\mfrak{p}}) \cong
C_{\mfrak{p}} \otimes_C \opn{Ext}_{C}^{i_0 - j_1}(N, M) . \] 
\end{proof}

\begin{rem}
Assume $A$ is noetherian. 
Proposition \ref{prop2.4} can be extended by replacing conditions 
(iii.a) and (iii.b) respectively with: 
(iii.a') $\mrm{H} B$ is a bounded essentially finite type 
$A$-algebra, and $L$ has finite projective dimension over $B$;
and (iii.b') $\mrm{H} C$ is a bounded essentially finite type 
$A$-algebra, $\mrm{H} N$ is a finitely generated $\mrm{H} C$ 
-module, and $N$ has finite flat dimension over $B$. 
The trick for (iii.a') is to localize on $\opn{Spec} \mrm{H}^0 B$
and to look at minimal semi-free resolutions of $L$. 
This trick also shows that 
$\opn{flat{.}dim}_B L = \opn{proj{.}dim}_B L$.
Details will appear elsewhere.
\end{rem}

%\newpage
%% ** section 2 **
\section{The Squaring Operation}
\label{sec2}

In this section we introduce a key technical notion used in the
definition of rigidity, namely the squaring operation. This 
operation is easy to define when the base ring is a field 
(see Corollary \ref{cor2.8}), but otherwise there are
complications, due to torsion. 
We solve the problem using DG algebras.

Recall that for a DG algebra $A$ the derived category of DG 
modules is denoted by $\til{\msf{D}}(\cat{DGMod} A)$. If $A$ is a 
ring then 
$\til{\msf{D}}(\cat{DGMod} A) = \msf{D}(\cat{Mod} A)$.

Let $A$ be a ring, and let
$M \in \msf{D}(\cat{Mod} A)$. As explained earlier the derived 
tensor product 
$M \otimes^{\mrm{L}}_{A} M \in \msf{D}(\cat{Mod} A)$
is defined to be
$M \otimes^{\mrm{L}}_{A} M := \til{M} \otimes_{A} \til{M}$,
where $\til{M} \to M$ is any K-flat resolution of $M$.
If $M \in \cat{DGMod} B$ for some DG $A$-algebra $B$, then we 
would like to be able to make
$M \otimes^{\mrm{L}}_{A} M$ into an object of 
$\til{\msf{D}}(\cat{DGMod} B \otimes_{A} B)$. 
But this is not always possible, at least not in any obvious way,
due to torsion. 
(For instance take $A := \mbb{Z}$ and $M = B := \mbb{Z} / (2)$). 
Fortunately there is a way to get around this problem.

Recall that a semi-free DG $A$-algebra $\til{B}$
is also semi-free as DG $A$-module, and hence it 
is a K-flat DG $A$-module. Thus 
Proposition \ref{prop1.2} implies that there exist K-flat DG 
algebra resolutions $A \to \til{B} \to B$ of $A \to B$.

\begin{lem} 
Let $A$ be a ring and $B$ a DG $A$-algebra. Suppose
$A \to \til{B} \to B$ is any K-flat DG 
algebra resolution of $A \to B$.
Then the \tup{(}non-additive\tup{)} functor
$\til{\msf{D}}(\cat{DGMod} B) \to \msf{D}(\cat{Mod} A)$,
$M \mapsto M \otimes^{\mrm{L}}_{A} M$,
factors canonically through \linebreak
$\til{\msf{D}}(\cat{DGMod} \til{B} \otimes_{A} \til{B})$. 
\end{lem}

\begin{proof}
Choose any quasi-isomorphism $\til{M} \to M$ in 
$\cat{DGMod} \til{B}$ with $\til{M}$ K-flat over $A$. This is 
possible since any K-flat DG $\til{B}$-module is 
K-flat over $A$. We get
$M \otimes^{\mrm{L}}_{A} M = \til{M} \otimes_{A} \til{M}
\in \til{\msf{D}}(\cat{DGMod} \til{B} \otimes_{A} \til{B})$. 
\end{proof}

\begin{thm} \label{thm1.1}
Let $A$ be a ring, let $B$ be a $A$-algebra 
and let $M$ be a DG $B$-module. 
Choose a K-flat DG algebra resolution 
$A \to \til{B} \to B$ of $A \to B$. 
Then the object 
\[ \opn{Sq}_{B / A} M :=
\mrm{RHom}_{\til{B} \otimes_{A} \til{B}}(B, 
M \otimes^{\mrm{L}}_{A} M)
\in \msf{D}(\cat{Mod} B) , \]
where the $B$-module structure is via the action on the first 
argument of $\opn{RHom}$, is independent of the choice of
resolution $A \to \til{B} \to B$.
\end{thm}

The functorial property of the assignment 
$M \mapsto  \opn{Sq}_{B / A} M$ is stated in 
Theorem \ref{thm1.2} below.

\begin{proof}
The idea for the proof was communicated to us by Bernhard Keller.
Choose some semi-free DG algebra 
resolution $A \to \til{B}' \to B$ of $A \to B$
(see Proposition \ref{prop1.2}).
We will show that there is a canonical isomorphism
\[ \mrm{RHom}_{\til{B} \otimes_{A} \til{B}}
(B, M \otimes^{\mrm{L}}_{A} M)
\iso 
\mrm{RHom}_{\til{B}' \otimes_{A} \til{B}'}
(B, M \otimes^{\mrm{L}}_{A} M) \]
in $\msf{D}(\cat{Mod} B)$. 

Let us choose a K-projective resolution $\til{M} \to M$ over
$\til{B}$, and a K-injective resolution
$\til{M} \otimes_{A} \til{M} \to \til{I}$
over $\til{B} \otimes_{A} \til{B}$. So
\[ \mrm{RHom}_{\til{B} \otimes_{A} \til{B}}(B, 
M \otimes^{\mrm{L}}_{A} M) =
\mrm{Hom}_{\til{B} \otimes_{A} \til{B}}
(B, \til{I}) . \]
Likewise let's choose resolutions 
$\til{M}' \to M$ and 
$\til{M}' \otimes_{A} \til{M}' \to \til{I}'$
over $\til{B}'$ and $\til{B}' \otimes_{A} \til{B}'$ 
respectively. 

By Proposition \ref{prop1.3} there is a DG algebra 
quasi-isomorphism $u_0 : \til{B}' \to \til{B}$ that's 
compatible with the quasi-isomorphisms to $B$. 
By the categorical properties of $K$-projective resolutions there
is an $\til{B}'$-linear quasi-isomorphism
$\phi_0 : \til{M}' \to \til{M}$, that's compatible  up
to homotopy with the quasi-isomorphisms to $M$. We obtain an 
$\til{B}' \otimes_{A} \til{B}'$ -linear quasi-isomorphism 
$\phi_0 \otimes \phi_0 : \til{M}' \otimes_{A} \til{M}'
\to \til{M} \otimes_{A} \til{M}$. 
Next by the categorical properties of $K$-injective resolutions 
there is an $\til{B}' \otimes_{A} \til{B}'$ -linear
quasi-isomorphism 
$\psi_0 : \til{I} \to \til{I}'$ 
that's compatible up to homotopy with the quasi-isomorphisms from 
$\til{M}' \otimes_{A} \til{M}'$. 
We thus get an $B$-linear homomorphism
\[ \chi_0 : \opn{Hom}_{\til{B} \otimes_{A} \til{B}} (B, \til{I}) 
\to \opn{Hom}_{\til{B}' \otimes_{A} \til{B}'}
(B, \til{I}') . \]
Proposition \ref{prop1.1} shows that $\chi_0$ is in fact an 
isomorphism in $\msf{D}(\cat{Mod} B)$. 

Now suppose $u_1 : \til{B}' \to \til{B}$,
$\phi_1 : \til{M}' \to \til{M}$ and
$\psi_1 : \til{I} \to \til{I}'$ 
are other choices of quasi-isomorphisms of the same respective 
types as $u_0$, $\phi_0$ and $\psi_0$. Then we get an induced 
isomorphism
\[ \chi_1 : \opn{Hom}_{\til{B} \otimes_{A} \til{B}} (B, \til{I}) 
\to \opn{Hom}_{\til{B}' \otimes_{A} \til{B}'}
(B, \til{I}') \]
in $\msf{D}(\cat{Mod} B)$. We shall prove that
$\chi_1 = \chi_0$. 

Here we have to introduce an auxiliary DG $A$-module 
$\mrm{C}(\til{M})$, the cylinder module. As graded module one has
\[ \mrm{C}(\til{M}) := 
\left[ \begin{matrix} \til{M} & \til{M}[-1] \\ 
0 & \til{M} \end{matrix} \right] , \] 
a triangular matrix module, and the differential is
\[ \mrm{d} \bigl( \left[ \begin{matrix} m_0 & n \\ 0 & m_1 
\end{matrix} \right] \bigr) :=
\left[ \begin{matrix} 
\mrm{d}(m_0) & \ m_0 - m_1 -\mrm{d}(n) \\ 
0 & \mrm{d}(m_1) \end{matrix} \right] \]
for $m_0, m_1, n \in \til{M}$.
There are DG module quasi-isomorphisms
$\epsilon : \til{M} \to \mrm{C}(\til{M})$
and
$\eta_0, \eta_1 : \mrm{C}(\til{M}) \to \til{M}$,
with formulas
\[ \epsilon(m) := 
\left[ \begin{matrix} m & 0 \\ 0 & m \end{matrix} \right]
\text{ and }
\eta_i \bigl( \left[ \begin{matrix} m_0 & n \\ 0 & m_1 
\end{matrix} 
\right] \bigr) := m_i . \] 
The cylinder module $\mrm{C}(\til{M})$ is a DG module over
$\til{B}$ by the formula
\[ a \cdot \left[ \begin{matrix} m_0 & n \\ 0 & m_1 
\end{matrix} \right] :=
\left[ \begin{matrix} a m_0 & a n \\ 0 & a m_1 
\end{matrix} \right] . \]

There is a quasi-isomorphism of DG $\til{B}$-modules
\[ \mrm{C}(\til{M}) \to 
\left[ \begin{matrix} \til{M} & M[-1] \\ 0 & \til{M} 
\end{matrix} \right] \]
which is the identity on the diagonal elements, and the given 
quasi-isomorphism $\til{M} \to M$ in the corner.
The two $\til{B}'$-linear quasi-isomorphisms $\phi_0$ and $\phi_1$ 
fit into an $\til{B}'$-linear quasi-isomorphism
\[ \til{M}' 
\xar{ \left[ \begin{smallmatrix} \phi_0 & 0 \\ 0 & \phi_1 
\end{smallmatrix} \right] } 
\left[ \begin{matrix} \til{M} & M[-1] \\ 
0 & \til{M} \end{matrix} \right] . \]
Since $\til{M}'$ is K-projective over $\til{B}'$ we can lift 
$\left[ \begin{smallmatrix} \phi_0 & 0 \\ 0 & \phi_1 
\end{smallmatrix} \right]$
to a quasi-isomorphism
$\phi : \til{M}' \to \mrm{C}(\til{M})$
such that $\eta_i \circ \phi = \phi_i$ up to homotopy.

Let's choose a K-injective resolution 
$\mrm{C}(\til{M}) \otimes_{A} \mrm{C}(\til{M}) \to \til{K}$ 
over $\til{B} \otimes_{A} \til{B}$. 
Then for $i = 0, 1$ we have a diagram
\[ \UseTips \xymatrix @C=5ex @R=5ex {
\til{I}' 
& \til{K}
\ar[l]_{\psi}
&  \til{I} 
\ar[l]_{\beta_i} \\
\til{M}' \otimes_{A} \til{M}'
\ar[u]
\ar[r]^(0.42){\phi \otimes \phi}
& \mrm{C}(\til{M}) \otimes_{A} \mrm{C}(\til{M})
\ar[u]
\ar[r]^(0.57){\eta_i \otimes \eta_i}
& \til{M} \otimes_{A} \til{M}
\ar[u] } \]
that's commutative up to homotopy. Here $\psi$ and $\beta_i$ are
some DG module homomorphisms, which 
exist due to the K-injectivity of $\til{I}'$ and 
$\til{K}$ respectively. Because 
$\phi_i \otimes \phi_i = (\eta_i \otimes \eta_i) \circ
(\phi \otimes \phi)$ up to homotopy,
and $\til{I}'$ is K-injective, it follows that the 
$\til{B}' \otimes_{A} \til{B}'$ -linear DG module 
quasi-isomorphisms $\psi \circ \beta_i$ and $\psi_i$
are homotopic. Therefore in order to prove that 
$\chi_0 = \chi_1$ it suffices to prove that the two isomorphisms
in $\msf{D}(\cat{Mod} B)$
\[ \theta_0, \theta_1 : 
\opn{Hom}_{\til{B} \otimes_{A} \til{B}}
(B, \til{I}) \to
\opn{Hom}_{\til{B} \otimes_{A} \til{B}} 
(B, \til{K}) , \]
that are induced by $\beta_0, \beta_1$ respectively, are equal.
 
For $i = 0, 1$ consider the diagram
\[ \UseTips \xymatrix @C=5ex @R=5ex {
\til{I} 
& \til{K}
\ar[l]_{\gamma}
&  \til{I} 
\ar[l]_{\beta_i} \\
\til{M} \otimes_{A} \til{M}
\ar[u]
\ar[r]^(0.42){\epsilon \otimes \epsilon}
& \mrm{C}(\til{M}) \otimes_{A} \mrm{C}(\til{M})
\ar[u]
\ar[r]^(0.57){\eta_i \otimes \eta_i}
& \til{M} \otimes_{A} \til{M}
\ar[u] } \] 
where $\gamma$ is some $\til{B} \otimes_{A} \til{B}$ -linear
DG module homomorphism, chosen so as to make the left square 
commute up to homotopy. As before, since 
$(\eta_i \otimes \eta_i) \circ (\epsilon \otimes \epsilon)
= \bsym{1}_{\til{M} \otimes_{A} \til{M}}$
it follows that 
$\gamma \circ \beta_i$ and $\bsym{1}_{\til{I}}$ 
are homotopic. Hence both $\theta_0$ and $\theta_1$ 
are inverses of the isomorphism
\[ \opn{Hom}_{\til{B} \otimes_{A} \til{B}}
(B, \til{K}) \iso 
\opn{Hom}_{\til{B} \otimes_{A} 
\til{B}}(B, \til{I}) \]
in $\msf{D}(\cat{Mod} B)$ induced by $\gamma$, so
$\theta_0 = \theta_1$.
\end{proof}

Suppose we are given a DG algebra homomorphism 
$B_1 \to B_2$; objects
$L_i, M_i \in \cat{DGMod} B_i$; and homomorphisms
$\psi : L_2 \to L_1$ and $\phi : M_1 \to M_2$ in 
$\cat{DGMod} B_1$. We then denote by
\[ (\psi, \phi) : \opn{Hom}_{B_1}(L_1, M_1) \to
\opn{Hom}_{B_2}(L_2, M_2) \]
the homomorphism 
$\chi \mapsto \phi \circ \chi \circ \psi$.
The same notation will be used for morphisms in derived 
categories.

\begin{thm} \label{thm1.2}
Let $A$ be a ring, let $B$ and $C$ be $A$-algebras, 
let $M \in \msf{D}(\cat{Mod} B)$ and 
$N \in \msf{D}(\cat{Mod} C)$.
Let $u : B \to C$ be a $A$-algebra homomorphism,
and let $\phi : N \to M$ be a morphism in $\msf{D}(\cat{Mod} B)$.
Suppose $A \to \til{B} \to B$ is a K-projective DG algebra 
resolution of $A \to B$, and $\til{B} \to \til{C} \to C$
is a K-projective DG algebra resolution  of $\til{B} \to C$.
Define a morphism 
\[ \opn{Sq}_{u / A}(\phi) : \opn{Sq}_{C / A} N \to
\opn{Sq}_{B / A} M \]
in $\msf{D}(\cat{Mod} B)$ by the formula
\[ \begin{aligned}
\opn{Sq}_{u / A}(\phi) :=
(u, \phi \otimes \phi) : \,
& \mrm{RHom}_{\til{C} \otimes_{A} \til{C}}(C, 
N \otimes^{\mrm{L}}_{A} N) \\
& \quad \to \mrm{RHom}_{\til{B} \otimes_{A} \til{B}}(B, 
M \otimes^{\mrm{L}}_{A} M) . 
\end{aligned} \]
Then the morphism $\opn{Sq}_{u / A}(\phi)$
is independent of the resolutions 
$A \to \til{B} \to B$ and $\til{B} \to \til{C} \to C$.
\end{thm}

The existence of these DG algebra resolutions is guaranteed by 
Propositions \ref{prop1.2} and \ref{prop1.3}.

\begin{proof}
Let's start by making the morphism 
$(u, \phi \otimes \phi)$ explicit. 
Note that $\til{C}$ is also K-projective relative 
to $A$, so it may be used to calculate $\opn{Sq}_{C / A} N$.
Let's choose DG module resolutions
$\til{M} \to M$, $\til{N} \to N$, 
$\til{M} \otimes_{A} \til{M} \to \til{I}$
and
$\til{N} \otimes_{A} \til{N} \to \til{J}$
by K-projective or K-injective DG modules over the appropriate DG 
algebras, as was done in the proof of Theorem \ref{thm1.1}.
Since $\til{N}$ is a K-projective DG $\til{B}$-module we get an 
actual DG module homomorphism 
$\til{\phi} : \til{N} \to \til{M}$
representing $\phi$. Therefore there is an 
$\til{B} \otimes_{A} \til{B}$ -linear DG module homomorphism
$\til{\phi} \otimes \til{\phi} : 
\til{N} \otimes_{A} \til{N} \to \til{M} \otimes_{A} 
\til{M}$.
Because $\til{I}$ is K-injective we obtain a DG module 
homomorphism $\psi : \til{J} \to \til{I}$
lifting 
$\til{N} \otimes_{A} \til{N} 
\xar{\til{\phi} \otimes \til{\phi}}
\til{M} \otimes_{A} \til{M} \to \til{I}$.
In this way we have obtained a homomorphism 
\[ (u, \psi) :
\mrm{Hom}_{\til{C} \otimes_{A} 
\til{C}}(C, \til{J}) \to
\mrm{Hom}_{\til{B} \otimes_{A} 
\til{B}}(B, \til{I}) \]
in $\cat{DGMod} B$, which represents $\opn{Sq}_{u / A}(\phi)$ once 
we pass to $\msf{D}(\cat{Mod} B)$.

Any other choices of K-projective and K-injective DG module
resolutions, and of the DG module homomorphisms 
$\til{\phi}$ and $\psi$, would yield a homomorphism that is 
homotopy equivalent to $(u, \psi)$. So the resulting morphism 
$\opn{Sq}_{u / A}(\phi)$ in 
$\msf{D}(\cat{Mod} B)$ will be unchanged.

It remains to prove $\opn{Sq}_{u / A}(\phi)$ is independent of 
the choice of DG algebra resolutions. 
Suppose we choose a
semi-free DG algebra resolution $A \to \til{B}' \to B$ of 
$A \to B$, and a  semi-free 
DG algebra resolution $\til{B}' \xar{\til{u}'} \til{C}' \to C$ of 
$\til{B}' \to C$. 
After choosing DG module resolutions 
$\til{M}' \to M$, $\til{N}' \to N$, 
$\til{M}' \otimes_{A} \til{M}' \to \til{I}'$
and
$\til{N}' \otimes_{A} \til{N}' \to \til{J}'$
by K-projective or K-injective DG modules over the appropriate DG 
algebras, we obtain a homomorphism
$\psi' : \til{J}' \to \til{I}'$ of DG 
$\til{B}' \otimes_{A} \til{B}'$ -modules, and morphism
\[ \opn{Sq}'_{u / A}(\phi) := (u, \psi') :
\mrm{Hom}_{\til{C}' \otimes_{A} \til{C}'}
(C, \til{J}') \to
\mrm{Hom}_{\til{B}' \otimes_{A} \til{B}'}
(B, \til{I}') \]
in $\cat{DGMod} B$. 

Applying Proposition \ref{prop1.3} twice we can find 
DG algebra homomorphisms $v_0$ and $w_0$ such that the 
diagram of DG algebra homomorphisms
\begin{equation} \label{eqn1.1}
\begin{array}{c}
\UseTips \xymatrix @C=5ex @R=4ex {
\til{B}'
\ar[r]^{v_0}
\ar[d]_{\til{u}'}
& \til{B}
\ar[d]_{\til{u}}
\ar[r] 
& B 
\ar[d]_{u} \\
\til{C}'
\ar[r]^{w_0}
& \til{C}
\ar[r] 
& C . } 
\end{array}
\end{equation}
is commutative. 
As in the proof of Theorem \ref{thm1.1} we pick 
quasi-isomorphisms 
$\psi_{M, 0} : \til{I} \to \til{I}'$ and
$\psi_{N, 0} : \til{J} \to \til{J}'$ over 
$\til{B}' \otimes_{A} \til{B}'$ and
$\til{C}' \otimes_{A} \til{C}'$ respectively. 
Then we get a commutative up to homotopy diagram 
\[ \UseTips \xymatrix @C=6ex @R=5ex {
\mrm{Hom}_{\til{B} \otimes_{A} \til{B}}
(B, \til{I}) 
\ar[r]^(0.47){\chi_{M, 0}} 
& \mrm{Hom}_{\til{B}' \otimes_{A} \til{B}'}
(B, \til{I}') \\
\mrm{Hom}_{\til{C} \otimes_{A} \til{C}}
(C, \til{J}) 
\ar[u]^{\opn{Sq}_{u / A}(\phi)}
\ar[r]^(0.47){\chi_{N, 0}} 
& \mrm{Hom}_{\til{C}' \otimes_{A} \til{C}'}
(C, \til{J}') 
\ar[u]_{\opn{Sq}'_{u / A}(\phi)}
} \] 
in $\cat{DGMod} B$, where the horizontal arrows are 
quasi-isomorphisms. If we were to choose another pair of DG 
algebra quasi-isomorphisms $v_1 : \til{B}' \to \til{B}$ and
$w_1 : \til{C}' \to \til{C}$ 
so as to make diagram \ref{eqn1.1} commutative, then according to
Theorem \ref{thm1.1} there would be equalities
$\chi_{M, 0} = \chi_{M, 1}$ and
$\chi_{N, 0} = \chi_{N, 1}$
of isomorphisms in $\msf{D}(\cat{Mod} B)$. Therefore 
$\opn{Sq}'_{u / A}(\phi) = \opn{Sq}_{u / A}(\phi)$
as morphisms in $\msf{D}(\cat{Mod} B)$.
\end{proof}

\begin{cor}
Let $A$ be a ring, let $B, C, D$ be $A$-algebras, 
let $M \in \msf{D}(\cat{Mod} B)$,  
$N \in \msf{D}(\cat{Mod} C)$ and
$P \in \msf{D}(\cat{Mod} D)$.
Let $B \xar{u} C \xar{v} D$ be $A$-algebra homomorphisms,
let $\phi : N \to M$ be a morphism in $\msf{D}(\cat{Mod} B)$, and
let $\psi : P \to N$ is a morphism in $\msf{D}(\cat{Mod} C)$.
\begin{enumerate}
\item There is equality
\[ \opn{Sq}_{v \circ u / A}(\phi \circ \psi) =
\opn{Sq}_{u / A}(\phi) \circ \opn{Sq}_{v / A}(\psi) \]
in 
$\opn{Hom}_{\msf{D}(\cat{Mod} B)}(\opn{Sq}_{D / A} P,
\opn{Sq}_{B / A} M)$.
\item If $B = C$, $M = N$, $u = \bsym{1}_B$ and 
$\phi = \bsym{1}_M$ \tup{(}the identity morphisms\tup{)}, then
$\opn{Sq}_{u / A}(\phi) = \bsym{1}_{\opn{Sq}_{B / A} M}$.
\end{enumerate}
\end{cor}

\begin{proof}
(1) Just choose a K-projective DG algebra resolution 
$\til{C} \to \til{D} \to D$ of
$\til{C} \to D$, and use the fact that
\[ \bigl( v \circ u, (\phi \circ \psi) \otimes (\phi \circ \psi)
\bigr) =
(u, \phi \otimes \phi) \circ
(v, \psi \otimes \psi) . \]

\medskip \noindent (2) Clear from the definitions.
\end{proof}

For the identity homomorphism $\bsym{1}_B : B \to B$
we write 
$\opn{Sq}_{B / A}(\phi) := \opn{Sq}_{\bsym{1}_B / A}(\phi)$.

\begin{dfn} \label{dfn2.8}
Let $A$ be a ring and let $B$ be an $A$-algebra.
The (nonlinear) functor
\[ \opn{Sq}_{B / A} : \msf{D}(\cat{Mod} B) \to
\msf{D}(\cat{Mod} B) \]
from Theorems \ref{thm1.1} and \ref{thm1.2}
is called the {\em squaring operation} over $B$ relative to $A$.
\end{dfn}

The next result explains the name ``squaring''. 

\begin{cor} \label{cor2.9}
In the situation of Theorem \tup{\ref{thm1.2}} let $c \in C$. Then 
\[ \opn{Sq}_{u / A}(c \phi) = c^2 \opn{Sq}_{u / A}(\phi) . \]
\end{cor}

\begin{proof}
It suffices to consider $u = \bsym{1}_C : C \to C$ and
$\phi = \bsym{1}_N : N \to N$. Choose any lifting of $c$ to 
$\til{c} \in \til{C}^0$. Then multiplication by
$\til{c} \otimes \til{c}$ on $\til{J}$ has the same effect on
$\mrm{Hom}_{\til{C} \otimes_{A} \til{C}}
(C, \til{J})$,
up to homotopy, as multiplication by $c^2$ on $C$.
\end{proof}

\begin{cor} \label{cor2.8} 
Suppose $B$ is a flat $A$-algebra, and $M$ is a bounded above 
complex of $B$-modules that are flat as $A$-modules. Then there 
is a functorial isomorphism
\[ \opn{Sq}_{B / A} M \cong \mrm{RHom}_{B \otimes_{A} B}
(B, M \otimes_{A} M) . \]
\end{cor}

\begin{proof}
This is because $B$ and $M$ are K-flat DG $A$-modules.
\end{proof}

\begin{rem}
One might be tempted to use the notation
$\mrm{RHom}_{B \otimes^{\mrm{L}}_{A} B}
(B, M \otimes_{A}^{\mrm{L}} M)$
instead of $\opn{Sq}_{B / A} M$. Indeed it is possible 
to make sense of the 
``DG algebra'' $B \otimes^{\mrm{L}}_{A} B$,
as an object of a suitable Quillen localization of the category of 
DG $A$-algebras. Cf.\ \cite{Hi} and \cite{Ta}, and also \cite{Qu}, where an 
analogous construction was made using simplicial algebras rather 
than DG algebras.  Then one should show that the 
triangulated category 
``$\til{\msf{D}}(\cat{DGMod} B \otimes^{\mrm{L}}_{A} B)$''
is well-defined, etc. See also \cite[Appendix V]{Dr}.
\end{rem}

%\newpage
%% ** section 3 **
\section{Essentially Smooth Homomorphisms} 
\label{sec3}

 From here on all rings are by default 
noetherian. We shall use notation such as $f^* : A \to B$ for a ring
homomorphism; so that $f : \opn{Spec} B \to \opn{Spec} A$ is the 
corresponding morphism of schemes. This will make our notation for 
various functors more uniform.
For instance restriction of scalars becomes 
$f_* : \cat{Mod} B \to \cat{Mod} A$, 
and extension of scalars (i.e.\ $M \mapsto B \otimes_A M$)
becomes
$f^* : \cat{Mod} A \to \cat{Mod} B$.
See also Definitions \ref{dfn3.1} and \ref{dfn2.6}.
Given another algebra 
homomorphism $g^* : B \to C$ we shall sometimes write 
$(f \circ g)^* := g^* \circ f^*$. 

In this section we present some results in commutative algebra. 
Recall that an $A$-algebra $B$ is called 
formally smooth (resp.\ formally
\'etale) if it has the lifting property (resp.\ the unique lifting 
property) for infinitesimal extensions. The $A$-algebra $B$ is called 
smooth (resp.\ \'etale) if it is finitely generated and
formally smooth (resp.\ formally \'etale). If $B$ is smooth over 
$A$ then it is flat, and $\Omega^1_{B / A}$ is a finitely generated
projective $B$-module. See 
\cite[Section $0{}_{\mrm{IV}}$.19.3]{EGA} 
and 
\cite[Section IV.17.3]{EGA} for details. 

\begin{dfn} \label{dfn3.3}
Let $A$ and $B$ be noetherian rings. A ring homomorphism 
$f^* : A \to B$ is called {\em essentially smooth} 
(resp.\ {\em essentially \'etale}) if it is of essentially finite 
type and formally smooth (resp.\ formally \'etale). 
In this case $B$ is called an essentially smooth
(resp.\ essentially \'etale) $A$-algebra.
\end{dfn}

Observe that smooth homomorphisms and localizations are 
essentially smooth.

\begin{prop} \label{prop3.2}
Let $f^* : A \to B$ be an essentially smooth homomorphism.
\begin{enumerate}
\item There is an open covering
$\opn{Spec} B = \bigcup_i \opn{Spec} B_i$ such that for every $i$ 
the homomorphism $A \to B_i$ is the composition of a smooth
homomorphism $A \to B_i^{\mrm{sm}}$ and a localization
$B_i^{\mrm{sm}} \to B_i$. 
\item $f^*$ is flat, and $\Omega^1_{B / A}$ is a 
finitely generated projective $B$-module.
\item $f^*$ is essentially \'etale if and only if 
$\Omega^1_{B / A} = 0$.
\item Let $g^* : B \to C$ be another essentially smooth 
homomorphism. Then 
$g^* \circ f^* : A \to C$ is also essentially smooth.
\end{enumerate}
\end{prop}

\begin{proof}
(1) Choose a finitely generated $A$-subalgebra 
$B^{\mrm{f}} \subset B$ such that $B$ is a localization of 
$B^{\mrm{f}}$. We can identify $U := \opn{Spec} B$ with a subset of 
$U^{\mrm{f}} := \opn{Spec} B^{\mrm{f}}$. Take a point $x \in U$, 
and let $y := f(x) \in \opn{Spec} A$. Then the local ring 
$\mcal{O}_{U^{\mrm{f}}, x} = \mcal{O}_{U, x} = B_x$ 
is a formally smooth $A_y$-algebra. According to 
\cite[Chapitre IV Th\'eor\`eme 17.5.1]{EGA}
there is an open neighborhood $W$ of $x$ in $U^{\mrm{f}}$ 
which is smooth over $\opn{Spec} A$. Choose 
an element $b \in B^{\mrm{f}}$ such 
that the localization $B^{\mrm{f}}[b^{-1}]$ 
satisfies
$x \in \opn{Spec} B^{\mrm{f}}[b^{-1}] \subset W$. 
Then $B^{\mrm{f}}[b^{-1}]$ is a smooth $A$-algebra,
$B[b^{-1}]$ is a localization of $B^{\mrm{f}}[b^{-1}]$,
$\opn{Spec} B[b^{-1}]$ is open in $\opn{Spec} B$, and
$x \in \opn{Spec} B[b^{-1}]$.
Finally let $i$ be an index corresponding to the point $x$, and
define $B_i^{\mrm{sm}} := B^{\mrm{f}}[b^{-1}]$ and 
$B_i := B[b^{-1}]$.

\medskip \noindent 
(2) This follows from (1). 

\medskip \noindent
(3)  See \cite[Chapitre $\mrm{0}_{\mrm{IV}}$ 
Proposition 20.7.4]{EGA}.

\medskip \noindent
(4) Both conditions in Definition \ref{dfn3.3} are transitive.  
\end{proof}

\begin{dfn}
Let $f^* : A \to B$ be an essentially smooth homomorphism.
If  $\opn{rank}_B \Omega^1_{B / A} = n$ 
then $f^*$ is called an
{\em essentially smooth homomorphism of relative dimension $n$}, 
and $B$ is called an {\em essentially smooth $A$-algebra 
of relative dimension $n$}.
\end{dfn}

By Proposition \ref{prop3.2}(3), 
an essentially \'etale homomorphism is the same as an 
essentially smooth homomorphism of relative dimension $0$.

\begin{prop} \label{prop3.4}
Suppose $f^* : A \to B$ and $g^* : B \to C$ are essentially 
smooth homomorphism of relative dimensions $m$ and $n$ 
respectively. Then $g^* \circ f^* : A \to C$ is an essentially 
smooth homomorphism of relative dimension $m + n$, and there is a 
canonical isomorphism of $C$-modules
$\Omega^{m+n}_{C / A} \cong \Omega^m_{B / A} \otimes_B
\Omega^n_{C / B}$.
\end{prop}

\begin{proof}
By \cite[Chapitre $\mrm{0}_{\mrm{IV}}$ Th\'eor\`eme 20.5.7]{EGA}
the sequence of $C$-modules
\[ 0 \to C \otimes_B \Omega^1_{B / A} \to
\Omega^1_{C / A} \to \Omega^1_{C / B} \to 0 \]
is split-exact. Choose any splitting
$s : \Omega^1_{C / B} \to \Omega^1_{C / A}$.
Then the homomorphism
\[ \Omega^m_{B / A} \otimes_B \Omega^n_{C / B} \to
\Omega^{m+n}_{C / A} , \]
\[ (\beta_1 \wedge \cdots \wedge \beta_m) \otimes
(\gamma_1 \wedge \cdots \wedge \gamma_n) \mapsto
\beta_1 \wedge \cdots \wedge \beta_m \wedge s(\gamma_1) 
\wedge \cdots \wedge s(\gamma_n) \]
is bijective and independent of the lifting $s$.
\end{proof}

Let us pause to review some facts about Koszul complexes. Suppose
$\bsym{a} = (a_1, \ldots, a_n)$ is a sequence 
of elements in a ring $A$. There is an associated
Koszul complex $\mbf{K}(A, \bsym{a})$, which is 
a super-commutative non-positive DG algebra. As graded algebra 
$\mbf{K}(A, \bsym{a})$ is the exterior algebra over the 
ring $A$ of the free module $\boplus_{i=1}^n A \til{a}_i$,
where $\til{a}_1, \ldots, \til{a}_n$ are variables of degree $-1$. 
(In the language of Section 1 we can say that 
$\mbf{K}(A, \bsym{a}) = A[\til{a}_1, \ldots, \til{a}_n]$,
the free super-commutative algebra on these odd generators.)
The differential is 
$\d(\til{a}_i) := a_i \in \mbf{K}^0(A, \bsym{a}) = A$. 

Suppose $\bsym{a}$ is a regular sequence in $A$. Let $J$ be 
the ideal generated by $\bsym{a}$, and let
$B := A / J$, the quotient ring. Then the augmentation
$\mbf{K}(A, \bsym{a}) \to B$ is a quasi-isomorphism, and so 
$\mbf{K}(A, \bsym{a})$ is a free resolution of $B$ as $A$-module.
It follows that for any $A$-module $M$ we have
\[  \opn{Ext}^p_{A}(B, M) =
\mrm{H}^p \opn{Hom}_{A} \bigl( \mbf{K}(A, \bsym{a}), M \bigr) . \]
Given an element $\mu \in M$ we define the generalized fraction
\[ \gfrac{\mu}{\bsym{a}} \in
\mrm{H}^n \opn{Hom}_{A} \bigl( \mbf{K}(A, \bsym{a}), M \bigr)
= \opn{Ext}^n_{A} (B, M ) \]
to be the cohomology class of the homomorphism
$\mbf{K}^{-n}(A, \bsym{a}) \to M$,
$\til{a}_1 \wedge \cdots \wedge \til{a}_n \mapsto \mu$.

The $B$-module $J / J^2$ is free, with basis the images
$\bar{a}_1, \ldots, \bar{a}_n$ of the regular sequence. 
Let us write
$\opn{det}(\bar{\bsym{a}}) := 
\bar{a}_1  \wedge \cdots \wedge \bar{a}_n
\in \bwedge^n_B (J / J^2)$.
Define the $B$-module
\[ \bsym{\omega}_{B / A} := 
\opn{Hom}_B \bigl(\bwedge^n_B (J / J^2) , B \bigr) , \]
which is free of rank $1$ with basis
$\frac{1}{\opn{det}(\bar{\bsym{a}})}$.
Then according to \cite[Proposition III.7.2]{RD} the map
\begin{equation} \label{eqn3.5}
\opn{Ext}^n_{A}(B, M) \to \bsym{\omega}_{B / A} \otimes_A M , 
\quad \gfrac{\mu}{\bsym{a}} \mapsto 
\frac{1}{\opn{det}(\bar{\bsym{a}})} \otimes \mu,
\end{equation}
is an isomorphism, and moreover it is independent of the regular 
sequence. Let us recall the argument: if $\bsym{a}'$ is another 
regular sequence generating $J$, then there is an invertible matrix
$\bsym{g} = [g_{i,j}]$ with entries in $A$
relating the two sequences, i.e.\ 
$a'_i = \sum_j g_{i,j} a_j$. Then
$\gfrac{\mu}{\bsym{a}'} = \opn{det}(\bsym{g})^{-1} 
\gfrac{\mu}{\bsym{a}}$,
and likewise  
$\frac{1}{\opn{det}(\bar{\bsym{a}}')} =
\opn{det}(\bsym{g})^{-1} \frac{1}{\opn{det}(\bar{\bsym{a}})}$.
Therefore the map is the same.

Clearly 
$\opn{Ext}^p_{A}(B, M) = 0$ for all $p > n$. Furthermore,
it is shown in \cite[Proposition III.7.2]{RD} 
that if $M$ is a flat $A$-module then 
$\opn{Ext}^p_{A}(B, M) = 0$ for all $p < n$.

 From here until Theorem \ref{thm3.1} (inclusive) we consider a 
flat ring homomorphism $A \to B$, and we let 
$B^{\mrm{e}} := B \otimes_A B$.
The diagonal embedding is
$\Delta : \opn{Spec} B \inj \opn{Spec} B^{\mrm{e}}$,
and we define
$J := \opn{Ker} (\Delta^* : B^{\mrm{e}} \to B)$.
Given an element $s \in B^{\mrm{e}}$ we write
$B[s^{-1}] := B[\Delta^*(s)^{-1}]$ and
\[ J [s^{-1}] := \opn{Ker} \bigl( B^{\mrm{e}}[s^{-1}] \to
B[s^{-1}] \bigr) \cong B^{\mrm{e}}[s^{-1}]
\otimes_{B^{\mrm{e}}} J . \]

\begin{lem} \label{lem3.2}
Assume $A \to B$ is essentially smooth of relative dimension $n$.
Then $\Delta$ is a regular closed immersion. 
Namely there exist finitely many elements $s_i \in B^{\mrm{e}}$
such that:
\begin{enumerate}
\rmitem{i} 
$\opn{Spec} B = \bigcup_i \opn{Spec} B[s_i^{-1}]$.
\rmitem{ii} For any $i$ the ideal $J[s_i^{-1}]$
is generated by a regular sequence of length $n$
in $B^{\mrm{e}}[s_i^{-1}]$.
\end{enumerate}
\end{lem}

\begin{proof}
Choose any $x \in \opn{Spec} B$. By Proposition 
\ref{prop3.2}(1) there are rings $B'$ and $B^{\mrm{sm}}$ such that
$A \to B^{\mrm{sm}}$ is smooth of relative dimension $n$; 
$B^{\mrm{sm}} \to B'$ is a 
localization; $\opn{Spec} B' \subset \opn{Spec} B$ is open; and
$x \in \opn{Spec} B'$. We may assume that 
$B' = B[r^{-1}]$ for some $r \in B$.

Since $A \to B^{\mrm{sm}}$ is smooth of relative dimension $n$,
according to \cite[Chapitre IV Proposition 17.12.4]{EGA}
there exists an element 
$t \in B^{\mrm{sm}} \otimes_A B^{\mrm{sm}}$
such that 
$x \in 
\opn{Spec} B^{\mrm{sm}}[t^{-1}]$,
and the ideal
\[ J^{\mrm{sm}}[t^{-1}] := 
\opn{Ker} \bigl( (B^{\mrm{sm}} \otimes_A B^{\mrm{sm}})[t^{-1}] \to
B^{\mrm{sm}}[t^{-1}] \bigr) \]
is generated by a regular sequence of length $n$. 

Now there is a ring homomorphism
\[ B^{\mrm{sm}} \otimes_A B^{\mrm{sm}} \to 
B[r^{-1}] \otimes_A B[r^{-1}] = 
B^{\mrm{e}}[(r \otimes r)^{-1}] . \]
Using it we can write 
$t = t_0 (r \otimes r)^{-l} \in B^{\mrm{e}}[(r \otimes r)^{-1}]$
for some non-negative integer $l$ and an element 
$t_0 \in B^{\mrm{e}}$.
Define $s_x := t_0 (r \otimes r) \in B^{\mrm{e}}$. Then the 
homomorphism 
$B^{\mrm{sm}} \otimes_A B^{\mrm{sm}} \to
B^{\mrm{e}}[s_x^{-1}]$
is flat, and the ideal
\[ J[s_x^{-1}]  \cong B^{\mrm{e}}[s_x^{-1}] 
\otimes_{B^{\mrm{sm}} \otimes_A B^{\mrm{sm}}} J^{\mrm{sm}}[t^{-1}] \]
is generated by the same regular sequence of length $n$. 

Going over all $x \in \opn{Spec} B$ we thus obtain a 
set $\{ s_x \}$ of elements of $B^{\mrm{e}}$ satisfying the 
two conditions (i,ii). Since $\opn{Spec} B$ is quasi-compact
we can select a finite subset. 
\end{proof}

\begin{rem}
In general the ring $B^{\mrm{e}}[s^{-1}]$ 
is not the same as the ring 
$B[s^{-1}] \otimes_A B[s^{-1}]$. 
For instance, 
take $A := \mbb{R}$ and $B := \mbb{C}$ (so $n = 0$ here).
Let $s := \mbf{i} \otimes 1 + 1 \otimes \mbf{i} \in B^{\mrm{e}}$.
Then $B^{\mrm{e}}[s^{-1}] = B[s^{-1}] = B$, 
but $B[s^{-1}] \otimes_A B[s^{-1}] = B^{\mrm{e}} \ncong B$.
\end{rem}

Take some $s \in B^{\mrm{e}}$, and let
$\bsym{b} = (b_1, \ldots, b_n)$ be a sequence of elements in 
$B^{\mrm{e}}[s^{-1}]$. we define
\[ \d(\bsym{b}) := \d(b_1) \wedge \cdots \wedge \d(b_n) \in
\Omega^n_{B^{\mrm{e}}[s^{-1}] / A} . \]

The de Rham complex $\Omega_{B/A}$ is a DG $A$-algebra (which 
lives in non-negative degrees!). The $A$-algebra homomorphism
$\mrm{p}_2^* : B \to B^{\mrm{e}}$, $b \mapsto 1 \otimes b$,
extends to a DG algebra homomorphism
$\mrm{p}_2^* : \Omega_{B/A} \to \Omega_{B^{\mrm{e}} / A}$.
Thus given an element $\beta \in \Omega^n_{B / A}$  
we obtain an element 
$\mrm{p}_2^*(\beta) \in \Omega^n_{B^{\mrm{e}} / A}$.
Clearly $(\Delta^* \circ \mrm{p}_2^*)(\beta) = \beta$.

\begin{lem} \label{lem3.4}
Assume $A \to B$ is essentially smooth of relative dimension $n$,
and $s \in B^{\mrm{e}}$ is such that the ideal
$J[s^{-1}]$
is generated by a regular sequence.
\begin{enumerate}
\item There is a unique 
$B[s^{-1}]$-linear isomorphism
\[ \Omega^{n}_{B[s^{-1}] / A} \iso
\opn{Ext}^n_{B^{\mrm{e}}[s^{-1}]} \bigl( B[s^{-1}], 
\Omega^{2n}_{B^{\mrm{e}}[s^{-1}] / A} \bigr) \]
such that
\[ \beta \mapsto 
\gfrac{\d(\bsym{b}) \wedge \mrm{p}_2^*(\beta)} 
{\bsym{b}} \]
for any regular sequence $\bsym{b} = (b_1, \ldots, b_n)$
generating $J[s^{-1}]$ and any $\beta \in \Omega^{n}_{B / A}$.
\item Let $M$ be any $B^{\mrm{e}}[s^{-1}]$-module. Then
\[ \opn{Ext}^p_{B^{\mrm{e}}[s^{-1}]} \bigl( B[s^{-1}], M \bigr)  
= 0 \]
for all $p > n$. Moreover, if $M$ is flat then this vanishing 
holds also for $p < n$.
\end{enumerate}
\end{lem}

\begin{proof}
(1) First we observe that there's an isomorphism of $B$-modules
\begin{equation} \label{eqn3.10}
\Omega^1_{B / A} \iso J / J^2, \quad
\d(b) \mapsto b \otimes 1 - 1 \otimes b \ \opn{mod} J^2
\end{equation}
for $b \in B$. In this way we get an isomorphism 
$\opn{Hom}_{B}(\Omega^{n}_{B / A}, B) \cong
\bsym{\omega}_{B / B^{\mrm{e}}}$.
The derivation 
$\d : B^{\mrm{e}} \to \Omega^{1}_{B^{\mrm{e}} / A}$
restricts to a function
$\d : J \to \Omega^{1}_{B^{\mrm{e}} / A}$
such that 
$\d(J^2) \subset J \cdot \Omega^{1}_{B^{\mrm{e}} / A}$, 
and this gives rise to is a split exact sequence of $B$-modules
\begin{equation} \label{eqn3.8}
0 \to J / J^2 \xar{\bar{\d}} 
B \otimes_{B^{\mrm{e}}} \Omega^{1}_{B^{\mrm{e}} / A}
\to \Omega^{1}_{B / A} \to  0 . 
\end{equation}
In this setup an element $b \in J \subset B^{\mrm{e}}$ 
will have images $\bar{b} \in J / J^2$ and
$\d(b) \in \Omega^{1}_{B^{\mrm{e}} / A}$; 
and these are related by the formula
$\bar{\d}(\bar{b}) = 1\otimes \d(b) \in 
B \otimes_{B^{\mrm{e}}} \Omega^{1}_{B^{\mrm{e}} / A}$.
Taking determinants (i.e.\ top degree exterior powers) in 
(\ref{eqn3.10}) and
(\ref{eqn3.8}) we obtain a canonical isomorphism 
\begin{equation} \label{eqn3.9}
B \otimes_{B^{\mrm{e}}} \Omega^{2n}_{B^{\mrm{e}} / A} \cong
\Omega^{n}_{B / A} \otimes_{B} \Omega^{n}_{B / A} .
\end{equation}
For any sequence $\bsym{b} = (b_1, \ldots, b_n)$ 
of elements in $J$ and differential 
form $\beta \in \Omega^{n}_{B / A}$ this isomorphism sends
$1 \otimes \bigl( \d(\bsym{b}) \wedge \mrm{p}_2^*(\beta) \bigr) 
\mapsto \opn{det}(\bar{\bsym{b}}) \otimes \beta$.

Now consider the isomorphism (\ref{eqn3.5}), 
but with updated entries. Using the isomorphism (\ref{eqn3.9})
it becomes
\[ \opn{Ext}^n_{B^{\mrm{e}}}(B, \Omega^{2n}_{B^{\mrm{e}} / A})
\cong \bsym{\omega}_{B / B^{\mrm{e}}} \otimes_{B^{\mrm{e}}}
\Omega^{2n}_{B^{\mrm{e}} / A} \cong 
\bsym{\omega}_{B / B^{\mrm{e}}} \otimes_{B}
\Omega^{n}_{B / A} \otimes_{B} \Omega^{n}_{B / A} . \]
After inverting $s$ we obtain a $B[s^{-1}]$-linear isomorphism
\begin{equation} \label{eqn3.6}
\begin{aligned}
& \opn{Ext}^n_{B^{\mrm{e}}[s^{-1}]}
(B[s^{-1}], \Omega^{2n}_{B^{\mrm{e}}[s^{-1}] / A}) \\
& \quad \cong \bsym{\omega}_{B[s^{-1}] / B^{\mrm{e}}[s^{-1}]} 
\otimes_{B[s^{-1}]}
\Omega^{n}_{B[s^{-1}] / A} \otimes_{B[s^{-1}]} 
\Omega^{n}_{B[s^{-1}] / A} . 
\end{aligned}
\end{equation}
Suppose the sequence $\bsym{b} = (b_1, \ldots, b_n)$ is 
regular and generates $J[s^{-1}]$. Then the form 
$\opn{det}(\bar{\bsym{b}})$ generates the $B[s^{-1}]$-module
\[ \bwedge^n_{B[s^{-1}]} \bigl( J[s^{-1}] / J[s^{-1}]^2 \bigr) 
= \bsym{\omega}_{B[s^{-1}] / B^{\mrm{e}}[s^{-1}]} . \]
Also any element of 
$B[s^{-1}] \otimes_{B^{\mrm{e}}[s^{-1}]}
\Omega^{2n}_{B^{\mrm{e}}[s^{-1}] / A}$ 
can be expressed as 
$\d(\bsym{b}) \wedge \mrm{p}_2^*(s^{-l} \beta)$ for some 
$\beta \in \Omega^{n}_{B / A}$ and $l \geq 0$. 
The isomorphism (\ref{eqn3.6}) sends
\[ \gfrac{ \d(\bsym{b}) \wedge \mrm{p}_2^*(s^{-l} \beta) }
{\bsym{b}} \mapsto
\frac{1}{\opn{det}(\bar{\bsym{b}})} \otimes
\opn{det}(\bar{\bsym{b}}) \otimes  s^{-l} \beta . \]
On the other hand the map
\[ s^{-l} \beta \mapsto \frac{1}{\opn{det}(\bar{\bsym{b}})} \otimes
\opn{det}(\bar{\bsym{b}}) \otimes  s^{-l} \beta  \]
is a $B[s^{-1}]$-linear isomorphism
\begin{equation} \label{eqn3.7}
\Omega^{n}_{B[s^{-1}] / A} \iso
\bsym{\omega}_{B[s^{-1}] / B^{\mrm{e}}[s^{-1}]} 
\otimes_{B[s^{-1}]}
\Omega^{n}_{B[s^{-1}] / A} \otimes_{B[s^{-1}]} 
\Omega^{n}_{B[s^{-1}] / A}
\end{equation}
which is evidently independent of the regular sequence $\bsym{b}$.
The isomorphism we want is the composition of (\ref{eqn3.6}) and
(\ref{eqn3.7}).

\medskip \noindent
(2) Take $s$ and $\bsym{b}$ as in part (1). 
We can use the Koszul complex
$\mbf{K} \bigl( B^{\mrm{e}}[s^{-1}], \bsym{b} \bigr)$
to calculate 
$\opn{Ext}^p_{B^{\mrm{e}}[s^{-1}]} \bigl( B[s^{-1}], M \bigr)$.
The assertions now follow from the general facts about Koszul 
complexes mentioned earlier; cf.\ \cite[Proposition III.7.2]{RD}.
\end{proof}

\begin{thm} \label{thm3.1}
Let $A \to B$ be an essentially smooth homomorphism of 
relative dimension $n$. Define 
$B^{\mrm{e}} := B \otimes_A B$. 
\begin{enumerate}
\item The $B^{\mrm{e}}$-module $B$ has finite projective 
dimension. 
\item One has
$\opn{Ext}^p_{B^{\mrm{e}}}(B, \Omega^{2n}_{B^{\mrm{e}} / A}) = 0$
for all $p \neq n$.
\item There is a unique $B$-linear isomorphism
\[ \Omega^n_{B / A} \iso 
\opn{Ext}^n_{B^{\mrm{e}}}(B, \Omega^{2n}_{B^{\mrm{e}} / A}) \]
which coincides with the isomorphism of Lemma \tup{\ref{lem3.4}(1)} 
after inverting any element $s \in B^{\mrm{e}}$ as described 
there.
\end{enumerate}
\end{thm}

\begin{proof}
(1, 2) We shall prove 
that for any $B^{\mrm{e}}$-module $M$ and any $p > n$ the module  
$\opn{Ext}^p_{B^{\mrm{e}}}(B, M)$ vanishes; and  if $M$ is flat 
there is also vanishing for $p < n$. Let 
$\opn{Spec} B = \bigcup_i \opn{Spec} B[s_i^{-1}]$ 
be the open covering from the Lemma \ref{lem3.2}. 
It suffices to show that 
$B[s_i^{-1}] \otimes_B \opn{Ext}^p_{B^{\mrm{e}}}(B, M) = 0$
for all $i$ and all $p$ in the corresponding range of integers.
But this was done in Lemma \ref{lem3.4}(2).

\medskip \noindent
(3) Uniqueness is clear, in view of Lemma \ref{lem3.2}. Regarding 
existence: due to the independence 
on regular sequences, the isomorphisms 
of Lemma \ref{lem3.4}(1) on the open sets $\opn{Spec} B[s_i^{-1}]$ 
can be glued.
\end{proof}

Finally a result about essentially \'etale homomorphisms. 

\begin{prop} \label{prop3.5}  
Let $A \to B$ be an essentially \'etale homomorphism,
and define $B^{\mrm{e}} := B \otimes_A B$.
Let $\mu : B^{\mrm{e}} \to B$ be the multiplication map, and let
$I := \opn{Ker}(\mu)$. 
We view $\opn{Hom}_{B^{\mrm{e}}}(B, B^{\mrm{e}})$
as an ideal of $B^{\mrm{e}}$, i.e.\ the annihilator of $I$.
Then there is a unique ring isomorphism
$\nu : B \times B' \iso B^{\mrm{e}}$
such that $\nu(B') = I$, 
$\nu(B) = \opn{Hom}_{B^{\mrm{e}}}(B, B^{\mrm{e}})$,
and $\mu \circ \nu : B \to B$ is the identity.
\end{prop}

\begin{proof}
Uniqueness: once we know that 
$B^{\mrm{e}} = I \oplus \opn{Hom}_{B^{\mrm{e}}}(B, B^{\mrm{e}})$
as $B^{\mrm{e}}$-modules, the ring $B'$ and the isomorphism $\nu$ 
are determined.

We need to construct $\nu$. First let's assume that $A \to B$ is 
\'etale (i.e.\ $B$ is a finitely generated $A$-algebra). 
Let $f^* : B \to B^{\mrm{e}}$ be the homomorphism
$f^*(b) := b \otimes 1$. Define $X := \opn{Spec} B^{\mrm{e}}$ and
$Y := \opn{Spec} B$; so $f : X \to Y$ is \'etale and separated.
The diagonal morphism $\Delta : Y \to X$ is a section of $f$, and
$\Delta^* = \mu$. 
According to \cite[Corollaire IV.17.9.3]{EGA} the morphism $\Delta$
is a closed an open immersion, so that 
$X = \Delta(Y) \coprod Z$ for some affine scheme
$Z = \opn{Spec} B'$. We get a ring isomorphism
$\nu :  B \times B' \iso B^{\mrm{e}}$ such that 
$\nu \circ \mu : B \to B$ is the identity. 
Let $\epsilon, \epsilon' \in B^{\mrm{e}}$ 
be the corresponding orthogonal idempotent 
elements, so that $\nu(B) = B^{\mrm{e}} \cdot \epsilon$ and 
$\nu(B') = B^{\mrm{e}} \cdot \epsilon'$.
Then 
$B^{\mrm{e}} \cdot \epsilon' = I$,
and the annihilator of $I$ in $B^{\mrm{e}}$ is precisely
$B^{\mrm{e}} \cdot \epsilon$.

Next lets look at a localization $B \to C$. Define
$C^{\mrm{e}} := C \otimes_A C$,
$\mu_C : C^{\mrm{e}} \to C$ the multiplication map, 
$I_C := \opn{Ker}(\mu_C)$ and
$C' :=  C^{\mrm{e}} \otimes_{B^{\mrm{e}}} B'$.
There is a factorization
$\nu_C : C \times C' \iso C^{\mrm{e}}$ 
gotten from $\nu$ by the tensor operation
$C^{\mrm{e}} \otimes_{B^{\mrm{e}}} -$.
Note that 
$I_C = C^{\mrm{e}} \cdot I$,
and so this decomposition can be characterized as
$\nu_C(C') = I_C$, $\nu_C(C)$ is the annihilator of $I_C$, and
$\mu_C \circ \nu_C : C \to C$ is the identity. 

Now for the general case. According to Proposition 
\ref{prop3.2}(1) we have an open covering
$\opn{Spec} B = \bigcup_i \opn{Spec} B_i$, such that for every $i$ 
the homomorphism $A \to B_i$ is the composition of an \'etale
homomorphism $A \to B_i^{\mrm{sm}}$ and a localization
$B_i^{\mrm{sm}} \to B_i$. As explained in the previous paragraph 
there are factorizations
$\nu_i : B_i \times B'_i \iso B_i \otimes_A B_i$,
and these are compatible with further localization. 
Therefore they can be glued into a global factorization
$\nu : B \times B' \iso B \otimes_A B$.
\end{proof}

%\newpage
%% ** section 4 **
\section{Rigid Complexes} 
\label{sec4}

In this section we introduce the main concept of the paper, namely 
{\em rigid complexes}. This concept is due to Van den Bergh 
\cite{VdB}. 

First a comment on bounded complexes. Let $B$ be a ring and
$M \in \msf{D}(\cat{Mod} B)$. If $M$ has bounded cohomology then,
after replacing it with the isomorphic complex 
$\tau^{\geq i} \tau^{\leq j} M$ for $i \gg 0$ and $j \ll 0$,
we can assume that $M \in \msf{D}^{\mrm{b}}(\cat{Mod} B)$.
Likewise for $\msf{D}^{+}(\cat{Mod} B)$.
Such considerations will be made implicitly throughout the paper.

For a noetherian ring $A$ we denote by 
$\msf{D}^{\mrm{b}}_{\mrm{f}}(\cat{Mod} A)$
the derived category of bounded complexes with finitely 
generated cohomologies. 

Let $A$ be a ring and let $B$ be an 
$A$-algebra. In Section \ref{sec2} we constructed a functor
$\opn{Sq}_{B / A} : \msf{D}(\cat{Mod} B) \to 
\msf{D}(\cat{Mod} B)$,
the squaring operation (see Definition \ref{dfn2.8}). 
When $A$ is a field one has the easy formula   
\[ \opn{Sq}_{B / A} M = \mrm{RHom}_{B \otimes_{A} B}
(B, M \otimes_{A} M)  \]
(see Corollary \ref{cor2.8}).
The squaring is functorial for algebra homomorphisms too. 
Given a homomorphism of algebras $f^* : B \to C$, 
complexes $M \in \msf{D}(\cat{Mod} B)$ and
$N \in \msf{D}(\cat{Mod} C)$, and a morphism
$\phi : N \to M$ in $\msf{D}(\cat{Mod} B)$, there is an induced 
morphism
$\opn{Sq}_{f^* / A}(\phi) : \opn{Sq}_{C / A} N \to 
\opn{Sq}_{B / A} M$
in $\msf{D}(\cat{Mod} B)$. 
Again when $A$ is a field the formula 
for $\opn{Sq}_{f^* / A}$ is obvious; complications arise only 
when the base ring $A$ is not a field. 

\begin{dfn} \label{dfn2.4}
Let $A$ be a ring, let $B$ be a noetherian
$A$-algebra, and let $M \in \msf{D}(\cat{Mod} B)$.
\begin{enumerate}
\item A {\em rigidifying isomorphism for $M$ relative to $A$}
is an isomorphism
\[ \rho : M \to \opn{Sq}_{B / A} M \]
in $\msf{D}(\cat{Mod} B)$. 
\item If $M \in \msf{D}^{\mrm{b}}_{\mrm{f}}(\cat{Mod} B)$ 
and it has finite flat dimension over $A$, then the pair 
$(M, \rho)$ in part \tup{(1)} is called a 
{\em rigid complex over $B$ relative to $A$}. 
\end{enumerate} 
\end{dfn}

\begin{exa}
Take $B = M := A$. Since $\opn{Sq}_{A / A} A = A$ it follows 
that $A$ has a tautological rigidifying isomorphism 
$\rho^{\mrm{tau}}_{A} : A \iso \opn{Sq}_{A / A} A$. 
We call $(A, \rho_{A}^{\mrm{tau}})$ the 
{\em tautological rigid complex over $A$ relative to $A$}. 
\end{exa}

\begin{dfn} \label{dfn2.5}
Let $A$ be a ring, let $B, C$ be noetherian $A$-algebras,
let $f^* : B \to C$ be a homomorphism of $A$-algebras, and 
let $(M, \rho_M)$ and $(N, \rho_N)$ be rigid complexes over $B$ 
and $C$ respectively, both relative to $A$. 
A morphism $\phi : N \to  M$ in 
$\msf{D}(\cat{Mod} B)$ is called a {\em rigid trace
morphism relative to $A$} if the diagram 
\[ \UseTips \xymatrix @C=5ex @R=5ex {
N 
\ar[r]^(0.4){\rho_N}
\ar[d]_{\phi}
& \opn{Sq}_{C / A} N
\ar[d]^{\opn{Sq}_{f^* / A}(\phi)} \\
M 
\ar[r]^(0.4){\rho_M}
& \opn{Sq}_{B / A} M
} \]
of morphisms in $\msf{D}(\cat{Mod} B)$ is commutative. 
If $B = C$ (and $f^*$ is the identity)
then we say $\phi : N \to M$ is a {\em rigid morphism over $B$
relative to $A$}.
\end{dfn}

It is easy to see that the composition of two rigid 
trace morphisms relative to $A$ is a again
rigid trace morphism 
relative to $A$. In particular, for 
a fixed $A$-algebra $B$ the rigid 
complexes over $B$ relative to $A$ form a category, which we 
denote by 
$\msf{D}^{\mrm{b}}_{\mrm{f}}(\cat{Mod} B)_{\mrm{rig} / A}$.

The importance of rigid complexes is captured by the next result.

\begin{thm} \label{thm4.1}
Let $A$ be a ring, let $B$ be a noetherian $A$-algebra, and let
$(M, \rho) \in 
\msf{D}^{\mrm{b}}_{\mrm{f}}(\cat{Mod} B)_{\mrm{rig} / A}$.
Assume the canonical homomorphism
$B \to  \opn{Hom}_{\msf{D}(\cat{Mod} B)}(M, M)$
is bijective.
Then the only automorphism of $(M, \rho)$ in 
$\msf{D}^{\mrm{b}}_{\mrm{f}}(\cat{Mod} B)_{\mrm{rig} / A}$
is the identity $\bsym{1}_M$.
\end{thm}

\begin{proof}
Let $\phi : M \to M$ be a rigid isomorphism. So
$\phi = b \bsym{1}_M$ for an invertible element $b \in B$. 
Then
\[ \opn{Sq}_{B / A}(\phi) \circ \rho = 
\rho \circ \phi =  \rho \circ (b \bsym{1}_M) = b \rho  \]
in 
$\opn{Hom}_{\msf{D}(\cat{Mod} B)}(M,
\opn{Sq}_{B / A} M)$. On the other hand, 
using Corollary \ref{cor2.9} we have
\[ \opn{Sq}_{B / A}(\phi) \circ \rho = 
\opn{Sq}_{B / A}(b \bsym{1}_M) \circ \rho 
= b^2 \opn{Sq}_{B / A}(\bsym{1}_M)  \circ \rho = 
b^2 \rho \circ \bsym{1}_M = b^2 \rho . \]
Hence $b = 1$.
\end{proof}

\begin{rem}
In our next paper \cite{YZ4} we prove the following result, which 
adds to the interest in rigid complexes. Suppose $A$ is a regular 
finite dimensional noetherian ring. Let $B$ be an essentially 
finite type $A$-algebra, and assume $\opn{Spec} B$ is connected 
and nonempty. Then, up to isomorphism, the category 
$\msf{D}^{\mrm{b}}_{\mrm{f}}(\cat{Mod} B)_{\mrm{rig} / A}$
contains {\em exactly two objects}: the zero complex, and a nonzero 
rigid complex $(R, \rho)$. The rigid complex $(R, \rho)$ 
satisfies the condition of Theorem \ref{thm4.1}, so it
has only one automorphism. Moreover the complex $R$ is 
{\em dualizing}.
\end{rem}

Let $A$ be a ring, let
$B \to C$ be a homomorphism of $A$-algebras, let
$M \in \msf{D}(\cat{Mod} B)$ and let
$N \in \msf{D}(\cat{Mod} C)$.
Choose a K-flat DG algebra resolution 
$A \to \til{B} \to B$ of $A \to B$, and then a 
K-flat DG algebra resolution $\til{B} \to \til{C} \to C$
of $\til{B} \to C$. This can be done by Proposition \ref{prop1.2}.
(If $A \to B$ and $B \to C$ are flat one may just 
take $\til{B} = B$ and $\til{C} = C$.)
Then there is a sequence of morphisms in 
$\msf{D}(\cat{Mod} C)$:
\begin{equation} \label{eqn4.4}
\begin{aligned}
& (\opn{Sq}_{B / A} M) \otimes^{\mrm{L}}_{B}
(\opn{Sq}_{C / B} N)  \\
& \quad = \mrm{RHom}_{\til{B} \otimes_{A} \til{B}}
(B, M \otimes^{\mrm{L}}_{A} M) 
\otimes^{\mrm{L}}_{\til{B}}
\mrm{RHom}_{\til{C} \otimes_{\til{B}} \til{C}}
(C, N \otimes^{\mrm{L}}_{\til{B}} N) \\
& \quad \to^{\diamondsuit}
\mrm{RHom}_{\til{C} \otimes_{\til{B}} \til{C}} \bigl(
C, N \otimes^{\mrm{L}}_{\til{B}} N \otimes^{\mrm{L}}_{\til{B}}
\mrm{RHom}_{\til{B} \otimes_{A} \til{B}}
(B, M \otimes^{\mrm{L}}_{A} M) \bigr) \\
& \quad \to^{\triangledown}
\mrm{RHom}_{\til{C} \otimes_{\til{B}} \til{C}} \Bigl(
C, \mrm{RHom}_{\til{B} \otimes_{A} \til{B}}
\bigl( B, (M \otimes^{\mrm{L}}_{\til{B}} N)
\otimes^{\mrm{L}}_{A}
(M \otimes^{\mrm{L}}_{\til{B}} N) \bigr) \Bigr) \\
& \quad \to^{\dag}
\mrm{RHom}_{\til{C} \otimes_{\til{B}} \til{C}} \Bigl(
C, \mrm{RHom}_{\til{C} \otimes_{A} \til{C}}
\bigl( \til{C} \otimes_{\til{B}} \til{C}, 
(M \otimes^{\mrm{L}}_{\til{B}} N)
\otimes^{\mrm{L}}_{A}
(M \otimes^{\mrm{L}}_{\til{B}} N) \bigr) \Bigr) \\
& \quad \to^{\heartsuit}
\mrm{RHom}_{\til{C} \otimes_{A} \til{C}} 
\bigl( C, 
(M \otimes^{\mrm{L}}_{\til{B}} N)
\otimes^{\mrm{L}}_{A}
(M \otimes^{\mrm{L}}_{\til{B}} N) \bigr) \\
& \quad = \opn{Sq}_{C / A}(M \otimes^{\mrm{L}}_{B} N)
 \end{aligned}
\end{equation}
defined as follows. The morphism $\diamondsuit$ is 
of the form
\begin{equation} \label{eqn4.5}
X \otimes^{\mrm{L}} \opn{RHom}(Y, Z) \to
\opn{RHom}(Y, X \otimes^{\mrm{L}} Z) .
\end{equation}
The morphism $\triangledown$ is a double application of 
(\ref{eqn4.5}).
The morphism $\dag$ is actually an isomorphism: it is 
Hom-tensor adjunction for the DG algebra homomorphism
$\til{B} \otimes_{A} \til{B} \to \til{C} \otimes_{A} \til{C}$,
plus the fact that 
\[ \til{C} \otimes_{\til{B}} \til{C} \cong
(\til{C} \otimes_{A} \til{C}) 
\otimes_{\til{B} \otimes_{A} \til{B}} B . \]
And the morphism $\heartsuit$ is also an isomorphism, being
Hom-tensor adjunction for the DG algebra homomorphism
$\til{C} \otimes_{A} \til{C} \to \til{C} 
\otimes_{\til{B}} \til{C}$.

\begin{lem} 
Let $A$ be a ring, let 
$f^* : B \to C$ be a homomorphism of $A$-algebras, let
$M \in \msf{D}(\cat{Mod} B)$ and let
$N \in \msf{D}(\cat{Mod} C)$. Then the morphism
\[ \smallsmile_{f^*; M, N} : 
(\opn{Sq}_{B / A} M) \otimes^{\mrm{L}}_{B} 
(\opn{Sq}_{C / B} N) \to
\opn{Sq}_{C / A} (M \otimes^{\mrm{L}}_{B} N) \]
in $\msf{D}(\cat{Mod} C)$ that was 
constructed in \tup{(\ref{eqn4.4})}
is independent of the resolutions 
$A \to \til{B} \to B$ and $\til{B} \to \til{C} \to C$.
\end{lem}

\begin{proof}
The homotopy arguments in the proofs of Theorem \ref{thm1.1}
apply here.
\end{proof}

The morphism $\smallsmile_{f; M, N}$ is called the 
{\em cup product}.  

\begin{lem} \label{lem4.1}  
The cup product $\smallsmile_{f; M, N}$ is functorial in $M$ and 
$N$. Namely suppose we are given a morphism
$\phi : M_1 \to M_2$ in $\msf{D}(\cat{Mod} B)$
and a morphism 
$\psi : N_1 \to N_2$ in $\msf{D}(\cat{Mod} C)$. Then the diagram
\begin{equation} \label{eqn4.7}
\UseTips \xymatrix @C=9ex @R=5ex {
(\opn{Sq}_{B / A} M_1) \otimes^{\mrm{L}}_{B} 
(\opn{Sq}_{C / B} N_1)
\ar[r]^(0.56){\smallsmile_{f; M_1, N_1}} 
\ar[d]_{\opn{Sq}_{B / A}(\phi) \ \otimes}
^{\opn{Sq}_{C / B}(\psi)}
& \opn{Sq}_{C / A} (M_1 \otimes^{\mrm{L}}_{B} N_1) 
\ar[d]^{\opn{Sq}_{C / A}(\phi \otimes \psi)} \\
(\opn{Sq}_{B / A} M_2) \otimes^{\mrm{L}}_{B} 
(\opn{Sq}_{C / B} N_2)
\ar[r]^(0.56){\smallsmile_{f; M_2, N_2}} 
& \opn{Sq}_{C / A} (M_2 \otimes^{\mrm{L}}_{B} N_2)
} 
\end{equation}
is commutative.
\end{lem}

\begin{proof}
Choose a semi-free DG algebra resolution 
$A \to \til{B} \to B$ of $A \to B$, and a semi-free DG 
algebra resolution $\til{B} \to \til{C} \to C$ of
$\til{B} \to C$. 
According to Theorem \ref{thm1.2}(3) we have representations
$\opn{Sq}_{B / A}(\phi) = (\bsym{1}_B, \phi \otimes \phi)$, 
$\opn{Sq}_{C / B}(\psi) = (\bsym{1}_C, \psi \otimes \psi)$ and
\[ \opn{Sq}_{C / A}(\phi \otimes \psi) = 
\bigl( (\bsym{1}_C, (\phi \otimes \psi) \otimes (\phi \otimes \psi)
\bigr) . \]
All the morphisms in (\ref{eqn4.4}) are compatible with 
$\phi$ and $\psi$. Hence the diagram \ref{eqn4.7} is 
commutative.
\end{proof}

\begin{thm} \label{thm2.5}
Let $A$ be a noetherian ring,
let $B$ and $C$ be essentially finite type $A$-algebras,
let $f^* : B \to C$ be an $A$-algebra homomorphism, 
$M \in \msf{D}^{\mrm{b}}_{\mrm{f}}(\cat{Mod} B)$
and 
$N \in \msf{D}^{\mrm{b}}_{\mrm{f}}(\cat{Mod} C)$.
Assume all three conditions \tup{(i)}, \tup{(ii)} and
\tup{(iii)} below hold.
\begin{enumerate}
\rmitem{i} The complex $M$ has finite flat dimension over $A$,
and $\mrm{H} \opn{Sq}_{B / A} M$ is \linebreak bounded. 
\rmitem{ii} The complex $N$ has finite flat dimension over $B$.
\rmitem{iii} Either \tup{(a)}, \tup{(b)} or \tup{(c)} is 
satisfied:
\begin{enumerate}
\rmitem{a} $B \to C$ is essentially smooth.
\rmitem{b} $M$ has finite flat dimension over $B$.
\rmitem{c} The the canonical morphism 
$B \to \mrm{RHom}_{B}(M, M)$ is an 
isomorphism.
\end{enumerate}
\end{enumerate}
Then:
\begin{enumerate}
\item The complex $M \otimes^{\mrm{L}}_{B} N$
is in $\msf{D}^{\mrm{b}}_{\mrm{f}}(\cat{Mod} C)$,
and it has finite flat dimension over $A$.
\item The cup product morphism
\[ \smallsmile_{f; M, N} :
(\opn{Sq}_{B / A} M) \otimes^{\mrm{L}}_{B} 
(\opn{Sq}_{C / B} N) \to
\opn{Sq}_{C / A} (M \otimes^{\mrm{L}}_{B} N) \]
is an isomorphism.
\end{enumerate}
\end{thm}

\begin{proof}
(1) The cohomology modules $\mrm{H}^i(M \otimes^{\mrm{L}}_{B} N)$
are finitely generated over $C$ because the rings are noetherian.
Since $M$ has finite flat dimension over $A$ and $N$ has 
finite flat dimension over $B$ (cf.\ Definition \ref{dfn2.4}), 
it follows that 
$M \otimes^{\mrm{L}}_{B} N \in \msf{D}(\cat{Mod} C)$ has 
finite flat dimension over $A$. In particular 
$\mrm{H}(M \otimes^{\mrm{L}}_{B} N)$ is bounded, and so 
we can assume that 
$M \otimes^{\mrm{L}}_{B} N \in 
\msf{D}^{\mrm{b}}_{\mrm{f}}(\cat{Mod} C)$.

\medskip \noindent
(2) As was already mentioned the morphisms $\dag$ and
$\heartsuit$ in (\ref{eqn4.4}) are automatically isomorphisms.

Let us consider the morphism $\triangledown$. In order to prove 
that the morphism
\begin{equation} \label{eqn4.6}
\begin{aligned}
& N \otimes^{\mrm{L}}_{\til{B}} N \otimes^{\mrm{L}}_{\til{B}} 
\mrm{RHom}_{\til{B} \otimes_{A} \til{B}} ( B, 
M \otimes^{\mrm{L}}_{A} M) \\
& \qquad \to \mrm{RHom}_{\til{B} \otimes_{A} \til{B}} \big( B, 
(M \otimes^{\mrm{L}}_{\til{B}} N) \otimes^{\mrm{L}}_{A} 
(M \otimes^{\mrm{L}}_{\til{B}} N) \big)  
\end{aligned}
\end{equation}
in $\til{\msf{D}}(\cat{DGMod} \til{C} \otimes_{A} \til{C})$
is an isomorphism, we can forget the 
$\til{C} \otimes_{A} \til{C}$ -module structure, and consider 
this as a morphism in $\msf{D}(\cat{Nod} A)$. 
According to Corollary \ref{cor1.1} the algebra 
$\mrm{H}^0(\til{B} \otimes_{A} \til{B}) \cong 
B \otimes_{A} B$ is noetherian, and each 
$\mrm{H}^i(\til{B} \otimes_{A} \til{B})$ is a finitely generated 
module over it. Since $N$ has finite flat dimension over 
$\til{B}$, and both $\mrm{H}(M \otimes^{\mrm{L}}_{A} M)$
and 
$\mrm{H} \bigl( (M \otimes^{\mrm{L}}_{\til{B}} N) 
\otimes^{\mrm{L}}_{A} M \bigr)$
are bounded, we can use Proposition \ref{prop2.4} twice, with its 
condition (iii.b), to conclude that (\ref{eqn4.6}) is an 
isomorphism.

Finally let's show that $\diamondsuit$ is an isomorphism.
The DG modules $N \otimes^{\mrm{L}}_{B} N$ and
$\opn{Sq}_{B / A} M$ have bounded cohomologies. 
If $B \to C$ is essentially smooth then 
$\til{C} \otimes_{\til{B}} \til{C} \to C \otimes_{B} C$ is a 
quasi-isomorphism, and moreover $C$ has finite projective 
dimension over $C \otimes_{B} C$. Thus under either condition 
(iii.a), (iii.b) or (iii.c) 
of the theorem we may apply Proposition \ref{prop2.4}, 
with its conditions (iii.a), (iii.b) or (iii.c) respectively, 
to deduce that $\diamondsuit$ is an isomorphism.
\end{proof}

The next result is an immediate consequence of the theorem.

\begin{cor} \label{cor4.2}
 Let $A$ be a noetherian ring,
let $B$ and $C$ be essentially finite type $A$-algebras,
let $f^* : B \to C$ be an $A$-algebra homomorphism, let
$(M, \rho_M) \in 
\msf{D}^{\mrm{b}}_{\mrm{f}}(\cat{Mod} B)_{\mrm{rig} / A}$
and let
$(N, \rho_N) \in 
\msf{D}^{\mrm{b}}_{\mrm{f}}(\cat{Mod} C)_{\mrm{rig} / B}$.
Assume condition \tup{(iii)} of Theorem \tup{\ref{thm2.5}} holds.
Then the complex 
$M \otimes^{\mrm{L}}_{B} N \in \msf{D}(\cat{Mod} C)$ 
has a unique rigidifying isomorphism
\[ \rho : M \otimes^{\mrm{L}}_{B} N \iso
\opn{Sq}_{C / A} (M \otimes^{\mrm{L}}_{B} N) , \]
such that the diagram
\[ \UseTips \xymatrix @C=9ex @R=5ex {
M \otimes^{\mrm{L}}_{B} N 
\ar[dr]^{\rho} 
\ar[d]_{\rho_M \otimes \rho_N}
\\
(\opn{Sq}_{B / A} M) \otimes^{\mrm{L}}_{B} 
(\opn{Sq}_{C / B} N)
\ar[r]_(0.55){\smallsmile_{f; M, N} }
& \opn{Sq}_{C / A} (M \otimes^{\mrm{L}}_{B} N)
} \]
is commutative.
\end{cor}

\begin{dfn} \label{dfn4.7}
In the situation of Corollary \ref{cor4.2}, the induced 
rigidifying isomorphism of $M \otimes^{\mrm{L}}_{B} N$
is denoted by $\rho_M \otimes \rho_N$. In addition we define
\[ (M, \rho_{M}) \otimes^{\mrm{L}}_{B} (N, \rho_{N}) :=
(M \otimes^{\mrm{L}}_{B} N, \rho_M \otimes \rho_N) 
\in \msf{D}^{\mrm{b}}_{\mrm{f}}(\cat{Mod} C)_{\mrm{rig} / A} . \]
\end{dfn}

\begin{cor} \label{cor4.1}
Let $f^* : B \to C$ be a homomorphism between 
essentially finite type $A$-algebras, let
$\phi : (M_1, \rho_{M_1}) \to (M_2, \rho_{M_2})$
be a morphism in 
$\msf{D}^{\mrm{b}}_{\mrm{f}}(\cat{Mod} B)_{\mrm{rig} / A}$,
and let 
$\psi : (N_1, \rho_{N_1}) \to (N_2, \rho_{N_2})$
be a morphism
$\msf{D}^{\mrm{b}}_{\mrm{f}}(\cat{Mod} C)_{\mrm{rig} / B}$.
Assume conditions \tup{(iii)} of Theorem \tup{\ref{thm2.5}} holds
with respect to the two pairs of complexes
$(M_1, N_1)$ and $(M_2, N_2)$, so that the rigidifying isomorphisms
$\rho_{M_1} \otimes \rho_{N_1}$ and $\rho_{M_2} \otimes \rho_{N_2}$
exist. Then the morphism
\[ \phi \otimes \psi : M_1 \otimes^{\mrm{L}}_{B} N_1 \to
M_2 \otimes^{\mrm{M}}_{B} N_2 \]
is a rigid morphism over $C$ relative to $A$.
\end{cor}

\begin{proof}
This is due to the functoriality of the cup product;
see Lemma \ref{lem4.1}.
\end{proof}

%\newpage
%% ** section 5 **
\section{Rigidity and Finite Homomorphisms} 
\label{sec.rig-finite}

In this section we introduce the inverse image operation for rigid 
complexes with respect to a finite ring homomorphism. 

\begin{dfn} \label{dfn3.1}
Let $f^* : A \to B$ be a ring homomorphism. 
\begin{enumerate}
\item Define a functor
$f^{\flat} : \msf{D}(\cat{Mod} A) \to \msf{D}(\cat{Mod} B)$
by $f^{\flat} M := \mrm{RHom}_{A}(B, M)$. 
\item Given $M \in \msf{D}(\cat{Mod} A)$ let 
$\opn{Tr}^{\flat}_{f; M} : f^{\flat} M \to M$
be the morphism $\phi \mapsto \phi(1)$.
This becomes a morphism of functors 
$\opn{Tr}^{\flat}_{f} : f_* f^{\flat} 
\to \bsym{1}_{\msf{D}(\cat{Mod} A)}$.
\end{enumerate}
\end{dfn}

Observe that if $g^* : B \to C$ is another ring homomorphism, then 
adjunction gives rise to a canonical isomorphism
\begin{equation} \label{eqn5.1}
g^{\flat} f^{\flat} M = 
\mrm{RHom}_{B} \bigl( C, \mrm{RHom}_{A}(B, M) \bigr)
\cong \mrm{RHom}_{A}(C, M) = (f \circ g)^{\flat} M . 
\end{equation}

\begin{thm} \label{thm2.2}
Let $A$ be a noetherian ring,
let $B$ and $C$ be essentially finite type $A$-algebras,
and let $f^* : B \to C$ be a finite $A$-algebra homomorphism.
Suppose we are given a rigid complex 
$(M, \rho) \in 
\msf{D}^{\mrm{b}}_{\mrm{f}}(\cat{Mod} B)_{\mrm{rig} / A}$,
such that $f^{\flat} M$ has finite flat dimension over $A$. 
\begin{enumerate}
\item The complex 
$f^{\flat} M \in \msf{D}^{\mrm{b}}_{\mrm{f}}(\cat{Mod} C)$ 
has an induced rigidifying isomorphism 
\[ f^{\flat}(\rho) : 
f^{\flat} M \iso \opn{Sq}_{C / A} f^{\flat} M . \]
The rigid complex
$f^{\flat} (M, \rho) := 
\bigl( f^{\flat} M, f^{\flat}(\rho) \bigr)$
depends functorially on $(M, \rho)$.
\item Suppose $g^* : C \to D$ is another finite homomorphism,
and $(f \circ g)^{\flat} M$ has finite flat dimension over $A$.
Then under the isomorphism
$(f \circ g)^{\flat} M \cong g^{\flat} f^{\flat} M$
of \tup{(\ref{eqn5.1})} one has
\[ g^{\flat}(f^{\flat}(\rho)) = (f \circ g)^{\flat} (\rho) . \]
\item The morphism 
$\opn{Tr}^{\flat}_{f; M} : f^{\flat} M \to M$
is a rigid trace morphism relative to $A$.
\end{enumerate}
\end{thm}

For the proof we will need a lemma.
The catch in this lemma is that the complex $P$ of flat 
$A$-module is bounded {\em below}, not above. 

\begin{lem} \label{lem2.4} 
Let $P$ and $N$ be bounded below complexes of $A$-modules. Assume 
that each $P^i$ is a flat $A$-module, and that $N$ has finite 
flat dimension over $A$. Then the canonical morphism 
$P \otimes^{\mrm{L}}_{A} N \to P \otimes_{A} N$
in $\msf{D}(\cat{Mod} A)$ is an isomorphism. 
\end{lem}

\begin{proof}
Choose a bounded flat resolution $Q \to N$ over $A$.
We have to show that 
$P \otimes_{A} Q \to P \otimes_{A} N$
is a quasi-isomorphism. 
Let $L$ be the cone of $Q \to N$.
It is enough to show that the complex $P \otimes_{A} L$
is acyclic. We note that $L$ is a bounded below acyclic 
complex and $P$ is a bounded below complex of flat modules. 
To prove that $\mrm{H}^i (P \otimes_{A} L) = 0$ for any given $i$ 
we might as well replace $P$ with a truncation 
$P' := ( \cdots \to P^{j_1 - 1} \to P^{j_1} \to 0 \to \cdots)$
for $j_1 \gg i$. Now $P'$ is K-flat, so 
$P' \otimes_{A} L$ is acyclic.
\end{proof}

\begin{proof}[Proof of the theorem]
(1) Let's pick a semi-free DG algebra resolution 
$A \to \til{B} \to B$ of $A \to B$. Next let's pick a 
K-projective DG algebra resolution $\til{B} \to \til{C} \to C$ of
$\til{B} \to C$, such that 
$\opn{und} \til{C} \cong \boplus_{i = -\infty}^0 
\opn{und} \til{B}[-i]^{\mu_i}$
with finite multiplicities $\mu_i$; see Proposition 
\ref{prop1.2}(3). Choose a bounded above semi-free resolution 
$P' \to M$ over $\til{B}$. Since $M$ has finite flat dimension over 
$A$ it follows that for $i \ll 0$ the truncated DG $\til{B}$-module 
$P := \tau^{\geq i} P'$ is a bounded complex of flat $A$-modules, 
and also $P \cong M$ in $\til{\msf{D}}(\cat{DGMod} \til{B})$. 

We have an isomorphism
$\mrm{Hom}_{\til{B}}(\til{C}, P) \cong \mrm{RHom}_{B}(C, M)
= f^{\flat} M$
in \linebreak $\til{\msf{D}}(\cat{DGMod} \til{C})$,
and an isomorphism
\[ \opn{Hom}_{\til{B} \otimes_{A} \til{B}}
(\til{C} \otimes_{A} \til{C}, P \otimes_{A} P) \cong
\opn{RHom}_{\til{B} \otimes_{A} \til{B}}
(\til{C} \otimes_{A} \til{C}, M \otimes^{\mrm{L}}_{A} M)  \]
in $\til{\msf{D}}(\cat{DGMod} \til{C} \otimes_{A} \til{C})$.
Because the multiplicities $\mu_i$ are finite and $P$ is bounded, 
the obvious DG module homomorphism
\[ \mrm{Hom}_{\til{B}}(\til{C}, P) \otimes_{A}
\mrm{Hom}_{\til{B}}(\til{C}, P) \to
\opn{Hom}_{\til{B} \otimes_{A} \til{B}}
(\til{C} \otimes_{A} \til{C}, P \otimes_{A} P) \]
is bijective. Now $\mrm{Hom}_{\til{B}}(\til{C}, P)$ is a bounded 
below complex of flat $A$-modules, which also has finite flat 
dimension over $A$. Therefore by Lemma \ref{lem2.4} we obtain
\[ \mrm{Hom}_{\til{B}}(\til{C}, P) \otimes_{A}
\mrm{Hom}_{\til{B}}(\til{C}, P) \cong 
\mrm{RHom}_{B}(C, M) \otimes^{\mrm{L}}_{A} \mrm{RHom}_{B}(C, M) 
\]
in $\til{\msf{D}}(\cat{DGMod} \til{C} \otimes_{A} \til{C})$.
We conclude that there is a functorial isomorphism 
\begin{equation} \label{eqn3.2} 
\mrm{RHom}_{B}(C, M) \otimes^{\mrm{L}}_{A} \mrm{RHom}_{B}(C, M)
\cong \opn{RHom}_{\til{B} \otimes_{A} \til{B}}
(\til{C} \otimes_{A} \til{C}, M \otimes^{\mrm{L}}_{A} M) 
\end{equation}
in $\til{\msf{D}}(\cat{DGMod} \til{C} \otimes_{A} \til{C})$.
(If $A$ is a field we may disregard 
the previous sentences, and just take $\til{B} := B$ and 
$\til{C} := C$.) We thus have a sequence of isomorphisms in 
$\msf{D}(\cat{Mod} C)$:
\begin{equation} \label{eqn2.1}
\begin{aligned}
\opn{Sq}_{C / A} f^{\flat} M 
& = \mrm{RHom}_{\til{C} \otimes_{A} \til{C}}
\bigl( C, \mrm{RHom}_B(C, M) \otimes^{\mrm{L}}_{A} 
\mrm{RHom}_B(C, M) \bigr) \\
& \cong^{\diamondsuit} \mrm{RHom}_{\til{C} \otimes_{A} \til{C}}
\bigl( C, \mrm{RHom}_{\til{B} \otimes_{A} \til{B}}
(\til{C} \otimes_{A} \til{C}, 
M \otimes^{\mrm{L}}_{A} M) \bigr) \\
& \cong^{\ddag} \mrm{RHom}_{\til{B} \otimes_{A} \til{B}}
( C, M \otimes^{\mrm{L}}_{A} M)  \\
& \cong^{\ddag} \mrm{RHom}_{B} \bigl( C, 
\mrm{RHom}_{\til{B} \otimes_{A} \til{B}}
( B, M \otimes^{\mrm{L}}_{A} M) \bigr)
= f^{\flat} \opn{Sq}_{B / A} M ,
\end{aligned}
\end{equation}
where the isomorphism marked $\diamondsuit$ is by (\ref{eqn3.2}), 
and the isomorphisms $\ddag$ come from the Hom-tensor
adjunction formula. The rigidifying isomorphism we want is the 
composition of 
$f^{\flat}(\rho) : f^{\flat} M \to f^{\flat} \opn{Sq}_{B / A} M$
with the isomorphism 
$\alpha_f : f^{\flat} \opn{Sq}_{B / A} M \iso 
\opn{Sq}_{C / A} f^{\flat} M$
we get from (\ref{eqn2.1}). 

\medskip \noindent (2)
Because of adjunction identities the diagram
\begin{equation} %\label{eqn2.5}
\UseTips \xymatrix @C=5ex @R=5ex {
(f \circ g)^{\flat} \opn{Sq}_{B / A} M
\ar[rr]^{\alpha_{f \circ g}}
\ar[d]_{\cong}
& & \opn{Sq}_{D / A} (f \circ g)^{\flat} M
\ar[d]_{\cong} \\
g^{\flat} f^{\flat} \opn{Sq}_{B / A} M
\ar[r]^{g^{\flat}(\alpha_f)}
& g^{\flat} \opn{Sq}_{C / A} f^{\flat} M
\ar[r]^{\alpha_g}
& \opn{Sq}_{D / A} g^{\flat} f^{\flat} M
} 
\end{equation}
in which the vertical arrows are the isomorphisms (\ref{eqn5.1}),
is commutative. Therefore 
$g^{\flat}(f^{\flat}(\rho)) = (f \circ g)^{\flat} (\rho)$.

\medskip \noindent (3)
Consider the diagram 
\begin{equation} \label{eqn2.5}
\UseTips \xymatrix @C=5ex @R=5ex {
f^{\flat} M
\ar[r]^(0.38){f^{\flat}(\rho)}
\ar[d]_{\phi}
& f^{\flat} \opn{Sq}_{B / A} M
\ar[d]_{\psi}
\ar[r]^{\alpha}
& \opn{Sq}_{C / A} f^{\flat} M
\ar[dl]!R^(0.5){\ \opn{Sq}_{f^* / A}(\phi)} \\
M 
\ar[r]^(0.4){\rho}
& \opn{Sq}_{B / A} M } 
\end{equation}
in which 
$\phi := \opn{Tr}^{\flat}_{f; M}$,
$N := \opn{Sq}_{B / A} M$ and
$\psi := \opn{Tr}^{\flat}_{f; N}$.
Since $\opn{Tr}^{}_{f} : f_* f^{\flat} \to \bsym{1}$
is a natural transformation, the square portion on diagram
(\ref{eqn2.5}) is commutative. 

As for the triangle portion of the diagram, according to Theorem 
\ref{thm1.2}(3) the morphism 
\[ \opn{Sq}_{f^* / A}(\phi) : \opn{Sq}_{C / A} f^{\flat} M \to
\opn{Sq}_{B / A} M \]
is just
\[ (f^*, \phi \otimes \phi) :
\mrm{RHom}_{\til{C} \otimes_{A} \til{C}}
\bigl( C, (f^{\flat} M) \otimes^{\mrm{L}}_{A} (f^{\flat} M)
\bigr) \to
\mrm{RHom}_{\til{B} \otimes_{A} \til{B}}
(B, M  \otimes^{\mrm{L}}_{A} M) . \]
Now going over the various isomorphisms in equation (\ref{eqn2.1})
we see that every one of them commutes with the obvious morphisms 
to 
$\mrm{RHom}_{\til{B} \otimes_{A} \til{B}}
(B, M  \otimes^{\mrm{L}}_{A} M)$.
Therefore the triangle portion of diagram
(\ref{eqn2.5}) is commutative, and so $\opn{Tr}^{\flat}_{f; M}$
is a rigid trace morphism.
\end{proof}

Suppose $B \to C$ is a ring homomorphism, 
$M \in \msf{D}(\cat{Mod} B)$ and $N \in \msf{D}(\cat{Mod} C)$.
Adjunction gives rise to a canonical isomorphism
\[ \mrm{RHom}_{B}(N, M) \cong 
\mrm{RHom}_{C} \bigl( N, \mrm{RHom}_{B}(C, M) \bigr) =
\mrm{RHom}_{C}(N, f^{\flat} M) . \] 
Given a morphism $\psi : N \to f^{\flat} M$ in
$\msf{D}(\cat{Mod} C)$, the corresponding morphism in 
$\msf{D}(\cat{Mod} B)$ is
$\opn{Tr}^{\flat}_{f; M} \circ\, \psi : N \to M$.

\begin{dfn} \label{dfn5.2}
Let $f^* : B \to C$ be a ring homomorphism, 
$M \in \msf{D}(\cat{Mod} B)$ 
and $N \in \msf{D}(\cat{Mod} C)$. A morphism
$\phi : N \to M$ in $\msf{D}(\cat{Mod} B)$ is called a 
{\em nondegenerate trace morphism} if the corresponding
morphism $N \to f^{\flat} M$ in 
$\msf{D}(\cat{Mod} C)$ is an isomorphism.
\end{dfn}

\begin{prop} \label{prop5.2}
Let $f^* : B \to C$ be a finite homomorphism between 
two essentially finite type $A$-algebras, let 
$(M, \rho_M) \in 
\msf{D}^{\mrm{b}}_{\mrm{f}}(\cat{Mod} B)_{\mrm{rig} / A}$
and
$(N, \rho_N) \in 
\msf{D}^{\mrm{b}}_{\mrm{f}}(\cat{Mod} C)_{\mrm{rig} / A}$.
Assume $N \cong f^{\flat} M$ and  
$\opn{Hom}_{\msf{D}(\cat{Mod} C)}(N, N) = C$. 
Then there exists a unique nondegenerate rigid trace morphism
$\phi : (N, \rho_N) \to (M, \rho_M)$ over $B$ relative to $A$.
\end{prop}

\begin{proof}
First note that in this case a morphism $N \to M$ is a 
nondegenerate trace morphism if and only if the corresponding 
morphism $N \to f^{\flat} M$ is a basis of the rank $1$ 
free $C$-module 
$\opn{Hom}_{\msf{D}(\cat{Mod} C)}(N, f^{\flat} M)$.

Because $N \cong f^{\flat} M$ this has finite flat dimension over 
$A$, so Theorem \ref{thm2.2}(1) applies, and the rigid complex
$f^{\flat}(M, \rho_M)$ exists.
Pick any isomorphism $\psi : N \iso f^{\flat} M$. Then
$\opn{Hom}_{\msf{D}(\cat{Mod} C)}(N, f^{\flat} M)$
is a free $C$-module with basis $\psi$. 
Since \linebreak
$\opn{Sq}_{C / A}(\psi) : \opn{Sq}_{C / A} N \to
\opn{Sq}_{C / A} f^{\flat} M$
is also an isomorphism, it follows that \linebreak
$\opn{Sq}_{C / A}(\psi) = u f^{\flat}(\rho_M) \circ \psi$
for some unique invertible $u \in C$.
Then
$u^{-1} \psi : (N, \rho_N) \to f^{\flat}(M, \rho_M)$
is a rigid isomorphism. By Theorem \ref{thm2.2}(2) the morphism
\[ \phi := \opn{Tr}^{\flat}_{f; M} \circ\, u^{-1} \psi : 
(N, \rho_N) \to (M, \rho_M) \]
is a nondegenerate rigid trace morphism.

Finally, since $u^{-1} \psi$ is the unique rigid isomorphism
$(N, \rho_N) \to f^{\flat}(M, \rho_M)$ it follows that
$\phi$ is the unique nondegenerate rigid trace morphism
$(N, \rho_N) \to (M, \rho_M)$.
\end{proof}

\begin{cor} \label{cor2.3}
In the situation of Theorem \tup{\ref{thm2.2}}, assume 
$\opn{Hom}_{\msf{D}(\cat{Mod} C)}(f^{\flat} M, f^{\flat} M) 
\linebreak = C$. 
Then $\opn{Tr}^{\flat}_{f; M}$ is the unique nondegenerate rigid 
trace morphism
$f^{\flat} (M, \rho_M) \to (M, \rho_M)$.
\end{cor}

\begin{proof}
Take $(N, \rho_N) := f^{\flat} (M, \rho_M)$
in Proposition \ref{prop5.2}.
\end{proof}

\begin{exa}
Suppose $f^* : A \to B$ is a finite flat ring homomorphism. 
Consider the rigid complex
$(N, \rho) := f^{\flat} (A, \rho_A^{\mrm{tau}}) \in 
\msf{D}^{\mrm{b}}_{\mrm{f}}(\cat{Mod} B)_{\mrm{rig} / A}$.
Since
$\opn{Sq}_{B / A} N \cong N$ it follows that
\[ \opn{Sq}_{B / A} N = 
\opn{Hom}_{B \otimes_A B}(B, N \otimes_A N) \subset
N \otimes_A N , \]
and in this way we may view $\rho$
as an $A$-linear homomorphism
$\rho : N \to N \otimes_A N$. 
It is not too hard to show (using Corollary \ref{cor2.3})
that $\rho$ makes $N$ into a 
coassociative cocommutative coalgebra over $A$, with counit 
$\opn{Tr}^{\flat}_{f; A} : N \to A$.
In fact, the coalgebra
$N$ is precisely the $A$-linear dual of the $A$-algebra $B$.
\end{exa}

%\newpage
%% ** section 6 **
\section{Rigidity and Essentially Smooth Homomorphisms} 
\label{sec.rig-smooth}

In this section we introduce the inverse image operation for rigid 
complexes with respect to an essentially smooth 
ring homomorphism. 

We shall need the following easy fact. 

\begin{lem} \label{lem3.1}
Suppose $B = \prod_{i = 1}^m B_i$, i.e.\ 
$\opn{Spec} B = \coprod_{i = 1}^m \opn{Spec} B_i$.
Then the functor $N \mapsto \prod_i (B_i \otimes_B N)$
is an equivalence 
$\msf{D}(\cat{Mod} B) \to 
\prod_{i} \msf{D}(\cat{Mod} B_i)$. 
\end{lem}

\begin{dfn} \label{dfn2.6}
Suppose $f^* : A \to B$ is an essentially smooth homomorphism 
between noetherian rings. Let
$\opn{Spec} B = \coprod_{i} \opn{Spec} B_i$
be the (finite) decomposition of $\opn{Spec} B$ into connected 
components. For each $i$ the $B_i$-module $\Omega^1_{B_i / A}$
is projective of constant rank, say $n_i$. Given 
$M \in \msf{D}(\cat{Mod} A)$ define
\[ f^{\sharp} M := \prod\nolimits_i 
(\Omega^{n_i}_{B_i / A}[n_i] \otimes_{A} M) . \]
This is a functor
$f^{\sharp} : \msf{D}(\cat{Mod} A) \to \msf{D}(\cat{Mod} B)$. 
\end{dfn}

Note that if $f^* : A \to B$ 
is essentially \'etale then one simply has
$f^{\sharp} M = B \otimes_A M$. 

\begin{thm} \label{thm2.4}
Let $A$ be a noetherian ring,
let $B$ and $C$ be essentially finite type $A$-algebras,
and let $f^* : B \to C$ be an essentially smooth $A$-algebra 
homomorphism. 
Let $(L, \rho) \in 
\msf{D}^{\mrm{b}}_{\mrm{f}}(\cat{Mod} B)_{\mrm{rig} / A}$. 
\begin{enumerate}
\item The complex $f^{\sharp} L$ has an induced rigidifying 
isomorphism 
\[ f^{\sharp}(\rho) : f^{\sharp} L \iso \opn{Sq}_{C / A} 
f^{\sharp} L . \]
We get a functor
$f^{\sharp} : 
\msf{D}^{\mrm{b}}_{\mrm{f}}(\cat{Mod} B)_{\mrm{rig} / A} \to
\msf{D}^{\mrm{b}}_{\mrm{f}}(\cat{Mod} C)_{\mrm{rig} / A}$. 
\item Let $(B, \rho_B^{\mrm{tau}})$ be the tautological rigid 
complex. Then under the standard isomorphism 
$f^{\sharp} L \cong L \otimes^{\mrm{L}}_{B} f^{\sharp} B$
one has 
$f^{\sharp}(\rho) = \rho \otimes f^{\sharp}(\rho_B^{\mrm{tau}})$.
\item Let $g^* : C \to D$ be an essentially smooth homomorphism.
Then under the isomorphism 
$(f \circ g)^{\sharp} L \cong g^{\sharp} f^{\sharp} L$
of Proposition \tup{\ref{prop3.4}} one has 
$(f \circ g)^{\sharp}(\rho) = g^{\sharp}  f^{\sharp} (\rho)$.
\end{enumerate}
\end{thm}

\begin{proof}
(1) In view of Lemma \ref{lem3.1} we might as well assume
$\Omega^1_{C / B}$ has constant rank $m$. 
There are ring homomorphisms
$\opn{p}_1^*, \opn{p}_2^* : C \to C \otimes_B C$,
namely $\opn{p}_1^*(c) := c \otimes 1$ and 
$\opn{p}_2^*(c) := 1 \otimes c$. These induce module homomorphisms
$\opn{p}_i^* : \Omega^{m}_{C / B} \to 
\Omega^{m}_{(C \otimes_B C) / B}$, 
and the resulting homomorphism
\[ \Omega^{m}_{C / B} \otimes_B \Omega^{m}_{C / B} \to
\Omega^{2m}_{(C \otimes_B C) / B}\ ,\
\gamma_1 \otimes \gamma_2 \mapsto \opn{p}_1^*(\gamma_1)
\wedge \opn{p}_2^*(\gamma_2)  \]
is bijective. So we can interpret Theorem \ref{thm3.1}(3) as a 
canonical rigidifying isomorphism for the complex 
$\Omega^{m}_{C / B}[m]$ relative to $B$, which we 
denote by $\rho_{\Omega}$. Thus we obtain an object
\[ (\Omega^{m}_{C / B}[m], \rho_{\Omega}) \in
\msf{D}^{\mrm{b}}_{\mrm{f}}(\cat{Mod} C)_{\mrm{rig} / B} . \]
Now using Theorem \ref{thm2.5}, under its condition (iii.a),
and Definition \ref{dfn4.7}, we can define the rigidifying 
isomorphism
$f^{\sharp}(\rho) := \rho \otimes \rho_{\Omega}$. 
Corollary \ref{cor4.1} guarantees functoriality.

\medskip \noindent (2)
This assertion is clear from part (1).

\medskip \noindent (3) 
We may assume $\Omega^1_{D / C}$ has constant rank $n$.
Let's introduce temporary notation
$C^{\mrm{e} / B} := C \otimes_B C$ and
$\mrm{E}(B, C) := 
\opn{Ext}^m_{C^{\mrm{e} / B}}
(C, \Omega^{2m}_{C^{\mrm{e} / B} / B})$,
so that
\[ \opn{Sq}_{C / B} \bigl( \Omega^{m}_{C / B}[m] \bigr) = 
\mrm{E}(B, C)[m] . \]
Likewise for $C \to D$ and $B \to D$.

Using part (2) of this theorem it suffices to prove that
\[ (f \circ g)^{\sharp}(\rho_B^{\mrm{tau}}) = 
f^{\sharp} (\rho_B^{\mrm{tau}}) \otimes
g^{\sharp} (\rho_C^{\mrm{tau}}) \]
as rigidifying isomorphisms for
$\Omega^{m+n}_{D / B}[m+n] = g^{\sharp}  f^{\sharp} B$ 
relative to $B$. In other words, we have to prove that the 
diagram of isomorphisms
\begin{equation} \label{eqn4.2}
\UseTips \xymatrix @C=5ex @R=5ex {
\Omega^{m}_{C / B} \otimes_C \Omega^{n}_{D / C}
\ar[r]^{} 
\ar[d]_{f^{\sharp} (\rho_B^{\mrm{tau}}) \otimes
g^{\sharp} (\rho_C^{\mrm{tau}})}
& \Omega^{m+n}_{D / B} 
\ar[d]^{(f \circ g)^{\sharp}(\rho_B^{\mrm{tau}})}\\
\mrm{E}(B, C) \otimes_C \mrm{E}(C, D)
\ar[r]
& \mrm{E}(B, D)
} 
\end{equation}
in which the top horizontal arrow is from Proposition 
\ref{prop3.4}, and the bottom horizontal arrow is  
the isomorphism of Theorem \ref{thm2.5}(2), shifted by $m+n$, is 
commutative. 

Let us define ideals
$J := \opn{Ker}(C^{\mrm{e} / B} \to C)$
and
$K := \opn{Ker}(D^{\mrm{e} / C} \to D)$.
Choose elements $s \in C^{\mrm{e} / B}$ 
and $t \in D^{\mrm{e} / C}$ such that the ideals
$J[s^{-1}]$ and $K[t^{-1}]$ are generated by regular sequences
$\bsym{c} = (c_1, \ldots, c_m)$ and
$\bsym{d} = (d_1, \ldots, d_n)$ respectively. 
Using the surjection $D^{\mrm{e} / B} \to D^{\mrm{e} / C}$
we lift $t$ to an element of $D^{\mrm{e} / B}$, which is also 
denoted by $t$. 
After possibly replacing $t$ by $st$ we can assume that 
$C^{\mrm{e} / B}[s^{-1}]$ goes into 
$D^{\mrm{e} / B}[t^{-1}]$.
In view of Lemma \ref{lem3.4} it suffices to verify that the 
diagram (\ref{eqn4.2}) commutes after inverting $t$. 

Let's write
\[ \mrm{E}(B, C)[s^{-1}] := 
\opn{Ext}^m_{C^{\mrm{e} / B}[s^{-1}]}
(C[s^{-1}], \Omega^{2m}_{C^{\mrm{e} / B}[s^{-1}] / B}) . \]
We can represent it using a Koszul complex:
\[ \mrm{E}(B, C)[s^{-1}]
 = \mrm{H}^m \opn{Hom}_{C^{\mrm{e} / B}[s^{-1}]} 
\Bigl( \mbf{K} \bigl( C^{\mrm{e} / B}[s^{-1}], \bsym{c} \bigr), 
\Omega^{2m}_{C^{\mrm{e} / B}[s^{-1}] / B} \Bigr) . \]
Likewise for $\mrm{E}(C, D)[t^{-1}]$ and
$\mrm{E}(B, D)[t^{-1}]$,
using the regular sequences $\bsym{d}$ and
$(\bsym{c}, \bsym{d})$ respectively.
Now let's return to the proof of Theorem \ref{thm2.5}(2) 
and track all the isomorphisms there.
But this time we have flat algebra homomorphisms
$B \to C \to D$ (instead of $A \to B \to C$), so there is no need 
to take DG algebra resolutions. On the other hand 
there are some localizations to remember. Thus, instead of 
$\til{B} \otimes_{A} \til{B}$ we now have 
$C^{\mrm{e} / B}[s^{-1}]$; instead of $B$ we now have
$C[s^{-1}]$; etc. And most importantly, the $\opn{RHom}$'s are
calculated using Koszul complexes, so elements of 
$\mrm{E}(B, C)[s^{-1}]$ etc.\ can be represented as generalized 
fractions. We see that the isomorphism
\[ \mrm{E}(B, C)[s^{-1}] \otimes_{C[s^{-1}]}
\mrm{E}(C, D)[t^{-1}] \iso \mrm{E}(B, D)[t^{-1}] \]
that we get from Theorem \ref{thm2.5}(2) sends
\[ \gfrac{\d(\bsym{c}) \wedge \mrm{p}_2^*(\gamma)}{\bsym{c}} 
\otimes 
\gfrac{\d(\bsym{d}) \wedge \mrm{p}_2^*(\delta)}{\bsym{d}} 
\mapsto
\gfrac{\d(\bsym{c}, \bsym{d}) \wedge 
\mrm{p}_2^*(\gamma \wedge \delta)}
{(\bsym{c}, \bsym{d})} . \]
Here $\gamma \in \Omega^{m}_{C / B}$ and
$\delta \in \Omega^{n}_{D / C}$ are arbitrary elements. 
Comparing this to Lemma \ref{lem3.4}(1) we conclude that indeed 
the diagram (\ref{eqn4.2}), localized, is commutative.
\end{proof}

\begin{dfn} \label{dfn4.6}
Let $A$ be a noetherian ring, and let
$f^* : A \to A'$ be an essentially \'etale ring 
homomorphism. For $M \in \msf{D}(\cat{Mod} A)$ let
$\mrm{q}^{\sharp}_{f; M} : M \to f^{\sharp} M = A' \otimes_A M$
be the morphism $m \mapsto 1 \otimes m$. As $M$ varies
this becomes a functorial morphism
$\mrm{q}^{\sharp}_{f} : \bsym{1}_{\msf{D}(\cat{Mod} A)} \to
f_* f^{\sharp}$.
\end{dfn}

In the situation of the definition above, given 
$M' \in \msf{D}(\cat{Mod} A')$, there is a canonical bijection
\[ \opn{Hom}_{\msf{D}(\cat{Mod} A)}(M, M') 
\cong \opn{Hom}_{\msf{D}(\cat{Mod} A')}(f^{\sharp} M, M') ,\
\phi \mapsto \bsym{1} \otimes \phi . \]
In particular, for 
$M' := f^{\sharp} M$, the morphism $\mrm{q}^{\sharp}_{f; M}$
corresponds to the identity $\bsym{1}_{M'}$.

\begin{dfn} \label{dfn6.1}
Let $A$ be a noetherian ring, 
let $A \to A'$ be an essentially \'etale ring homomorphism, let
$M \in \msf{D}(\cat{Mod} A)$ and
$M' \in \msf{D}(\cat{Mod} A')$.
A morphism $\phi : M \to M'$ in $\msf{D}(\cat{Mod} A)$
is called a {\em nondegenerate 
localization morphism} if the induced morphism
$\bsym{1} \otimes \phi : A' \otimes_A M \to M'$ is an isomorphism. 
\end{dfn}

\begin{dfn} \label{dfn2.3}
Let $A$ be a noetherian ring, 
let $B$ and $B'$ be essentially finite type $A$-algebras,
let $f^* : B \to B'$ be an essentially \'etale 
$A$-algebra homomorphism, let
$(M, \rho) \in 
\msf{D}^{\mrm{b}}_{\mrm{f}}(\cat{Mod} B)_{\mrm{rig} / A}$
and let
$(M', \rho') \in  
\msf{D}^{\mrm{b}}_{\mrm{f}}(\cat{Mod} B')_{\mrm{rig} / A}$.
A {\em rigid localization morphism} is a 
morphism $\phi : M \to M'$ in 
$\msf{D}(\cat{Mod} B)$,
such that the induced morphism
$\bsym{1} \otimes \phi : 
f^{\sharp} (M, \rho) \to (M', \rho')$ 
is a rigid morphism over $B'$ relative to $A$, in the sense of 
Definition \ref{dfn2.5}.
\end{dfn}

\begin{prop} \label{prop2.1}  
Let $A$ be a noetherian ring, 
let $B$ and $B'$ be essentially finite type $A$-algebras,
let $f^* : B \to B'$ be an essentially \'etale 
$A$-algebra homomorphism, and
let $(M, \rho) \in 
\msf{D}^{\mrm{b}}_{\mrm{f}}(\cat{Mod} B)_{\mrm{rig} / A}$. 
Define 
$(M', \rho') := f^{\sharp}(M, \rho)$.
Then:
\begin{enumerate}
\item The morphism $\mrm{q}^{\sharp}_{f; M} : M \to M'$ 
is a nondegenerate rigid localization morphism.
\item Moreover, if $\mrm{RHom}_{B}(M, M) = B$, then 
$\mrm{q}^{\sharp}_{f; M}$ is the unique 
nondegenerate rigid localization morphism $M \to M'$.
\end{enumerate}
\end{prop}

\begin{proof}
(1) Since the corresponding morphism $M' \to M'$
is the identity automorphism of $M'$, it is certainly rigid and 
nondegenerate.

\medskip \noindent
(2) Here we have
$\opn{Hom}_{\msf{D}(\cat{Mod} B')}(M', M') = B'$.  
The uniqueness of $\mrm{q}^{\sharp}_{f; M}$ is proved 
like in Corollary \ref{cor2.3}.
\end{proof}

\begin{thm}[Base Change] \label{thm6.3} 
Let $A$ be a noetherian ring, 
let $B, B'$ and $C$ be essentially finite type $A$-algebras,
and let $f^* : B \to C$ and $g^* : B \to B'$ be homomorphisms 
of $A$-algebras, with $g^*$ a localization. 
Define $C' := B' \otimes_B C$, so there are induced homomorphisms
${f'}^* : B' \to C'$ and $h^* : C \to C'$, with $h^*$ a 
localization; cf.\ diagram below.
Let 
$(M, \rho_M) \in 
\msf{D}^{\mrm{b}}_{\mrm{f}}(\cat{Mod} B)_{\mrm{rig} / A}$,
let
$(N, \rho_N) \in 
\msf{D}^{\mrm{b}}_{\mrm{f}}(\cat{Mod} C)_{\mrm{rig} / A}$,
and suppose 
\[ \phi : (N, \rho_N) \to (M, \rho_M) \] 
is a rigid trace morphism over $B$ relative to $A$. 
Define 
$M' := g^{\sharp} M$ and $N' := h^{\sharp} N$.
There is a morphism
$\phi' : N' \to M'$ in $\msf{D}(\cat{Mod} B')$
gotten by composing the canonical isomorphism
$N' = C' \otimes_C N \cong B' \otimes_B N = g^{\sharp} N$
with
$g^{\sharp}(\phi) : g^{\sharp} N \to g^{\sharp} M = M'$.
Then 
\[ \phi' : \bigl( N', h^{\sharp}(\rho_N) \bigr) \to
\bigl( M', g^{\sharp}(\rho_M) \bigr)  \]
is a rigid trace morphism over $B'$ relative to $A$.
\[ \UseTips \xymatrix @C=5ex @R=5ex {
B
\ar[r]^{f^*} 
\ar[d]_{g^*} 
& C
\ar[d]^{h^*} \\
B' 
\ar[r]_{{f'}^*}
& C'
} 
\qquad
\UseTips \xymatrix @C=5ex @R=5ex {
M 
\ar[d]_{\opn{q}^{\sharp}_{g; M}} 
& N
\ar[l]_{\phi}
\ar[d]^{\opn{q}^{\sharp}_{h; N}} \\
M' 
& N'
\ar[l]^{\phi'}
} \]
\end{thm}

Observe there is no particular assumption on $f^*$, and $\phi$ is 
not assumed to be nondegenerate. 
For the proof we shall need a few lemmas.

Choose a K-flat DG algebra resolution 
$A \to \til{B} \to B$ of $A \to B$, and 
a K-flat DG algebra resolution 
$\til{B} \to \til{B}' \to B'$ of $\til{B} \to B'$; 
see Proposition \ref{prop1.2}.
(If $A \to B$ and $B \to B'$ are flat then we may just take
$\til{B} := B$ and $\til{B}' := B'$.)
Define a morphism
\[ \lambda_{M} : \opn{Sq}_{B / A} M \to
\opn{Sq}_{B' / A} M' \]
in $\msf{D}(\cat{Mod} B)$ by composing the following morphisms:
\begin{equation} \label{eqn6.3}
\begin{aligned}
\opn{Sq}_{B / A} M 
& = \mrm{RHom}_{\til{B} \otimes_{A} \til{B}}(B, 
M \otimes^{\mrm{L}}_{A} M) \\
& \xar{(\bsym{1}, \mrm{q} \otimes \mrm{q})}
\mrm{RHom}_{\til{B} \otimes_{A} \til{B}}(B, 
M' \otimes^{\mrm{L}}_{A} M') \\
& \iso^{\diamondsuit}
\mrm{RHom}_{\til{B}' \otimes_{A} \til{B'}}(B',
M' \otimes^{\mrm{L}}_{A} M') \\
& = \opn{Sq}_{B' / A} M' .  
\end{aligned}
\end{equation}
In these formulas $\bsym{1}$ denotes the identity morphism
of $B$, and
$\mrm{q} := \opn{q}^{\sharp}_{g; M} : M \to M'$.
For the notation $(\bsym{1}, \mrm{q} \otimes \mrm{q})$ see 
text prior to Theorem \ref{thm1.2}. The isomorphism $\diamondsuit$ 
is by Hom-tensor adjunction with respect to the DG algebra 
homomorphism
$\til{B} \otimes_{A} \til{B} \to 
\til{B}' \otimes_{A} \til{B}'$,
plus the fact that
$B \otimes_{\til{B} \otimes_{A} \til{B}}
(\til{B}' \otimes_{A} \til{B}') \cong B'$.

\begin{lem} \label{lem6.3}
The morphism $\lambda_{M}$ 
is independent of the DG algebra 
resolutions $A \to \til{B} \to B$ and 
$\til{B} \to \til{B}' \to B'$.
\end{lem}

\begin{proof}
Choose a K-injective resolution
$M \otimes^{\mrm{L}}_{A} M \to I$
over $\til{B} \otimes_{A} \til{B}$, and a K-injective resolution
$M' \otimes^{\mrm{L}}_{A} M' \to I'$
over $\til{B}' \otimes_{A} \til{B}'$. Since 
$\til{B}' \otimes_{A} \til{B}'$ is K-flat over 
$\til{B} \otimes_{A} \til{B}$ it follows that $I'$ is a 
K-injective DG $\til{B} \otimes_{A} \til{B}$ -module. Therefore
$\opn{q} \otimes \opn{q} : M \otimes_{A} M \to 
M' \otimes_{A} M'$
extends to a homomorphism $\psi : I \to I'$ in 
$\cat{DGMod} \til{B} \otimes_{A} \til{B}$,
uniquely up to homotopy. The morphism $\lambda_{M}$ is 
represented by
\[ \mrm{Hom}_{\til{B} \otimes_{A} \til{B}}(B, I)
\xar{(\bsym{1}, \psi)}
\mrm{Hom}_{\til{B} \otimes_{A} \til{B}}(B, I') 
\iso^{\diamondsuit}
\mrm{Hom}_{\til{B}' \otimes_{A} \til{B'}}(B', I')  . \]
Now the homotopy argument from the proof of Theorem \ref{thm1.1} 
can be used.
\end{proof}

In the same way we define a morphism 
$\lambda_{N} : \opn{Sq}_{C / A} N \to
\opn{Sq}_{C' / A} N'$
in \linebreak $\msf{D}(\cat{Mod} C)$. 

\begin{lem} \label{lem6.4}
The diagram
\begin{equation} \label{eqn6.5}
\UseTips \xymatrix @C=7ex @R=5ex {
\opn{Sq}_{C / A} N
\ar[d]_{\opn{Sq}_{f^* / A}(\phi)}
\ar[r]^{\lambda_{N}}
& \opn{Sq}_{C' / A} N' 
\ar[d]^{\opn{Sq}_{{f'}^* / A}(\phi')} \\
\opn{Sq}_{B / A} M
\ar[r]^{\lambda_{M}} 
& \opn{Sq}_{B' / A} M'
}
\end{equation}
of morphisms in $\msf{D}(\cat{Mod} B)$ is commutative.
\end{lem}

\begin{proof}
Choose a semi-free DG algebra resolution 
$A \to \til{B} \to B$ of $A \to B$. Next choose semi-free
DG algebra resolutions
$\til{B} \xar{v} \til{B}' \to B'$ and 
$\til{B} \xar{u} \til{C} \to C$
of $\til{B} \to B'$ and $\til{B} \to C$ respectively. 
Define $\til{C}' := \til{B}' \otimes_{\til{B}} \til{C}$.
Because $B \to B'$ is flat we get a quasi-isomorphism
$\til{C}' \to C'$, and hence 
$\til{C} \xar{w} \til{C}' \to C'$
and 
$\til{B}' \xar{u'} \til{C}' \to C'$
are semi-free DG algebra resolutions of 
$\til{C} \to C'$ and $\til{B}' \to C'$ respectively. 
By definition of the morphisms $\lambda_{-}$ and $\opn{Sq}_{-}$
the diagram (\ref{eqn6.5}) is represented by
\[ \UseTips \xymatrix @C=9ex @R=5ex {
\mrm{RHom}_{\til{C} \otimes_{A} \til{C}}
(C, N \otimes^{\mrm{L}}_{A} N)
\ar[d]_{(f^*, \phi \otimes \phi)}
\ar[r]^(0.47){\opn{q} \otimes \opn{q}}
& \mrm{RHom}_{\til{C}' \otimes_{A} \til{C}'}
(C', N' \otimes^{\mrm{L}}_{A} N') 
\ar[d]^{({f'}^*, \phi' \otimes \phi')} \\
\mrm{RHom}_{\til{B} \otimes_{A} \til{B}}
(B, M \otimes^{\mrm{L}}_{A} M) 
\ar[r]^(0.47){\opn{q} \otimes \opn{q}}
& \mrm{RHom}_{\til{B}' \otimes_{A} \til{B}'}
(B', M' \otimes^{\mrm{L}}_{A} M')
} \]
(we are suppressing the adjunction isomorphisms $\diamondsuit$).
This latter diagram is trivially commutative. 
\end{proof}

\begin{lem} \label{lem6.5}
The diagram
\[ \UseTips \xymatrix @C=9ex @R=5ex {
M
\ar[d]_{\opn{q}} 
\ar[r]^(0.4){\rho_M}
& \opn{Sq}_{B / A} M 
\ar[d]^{\lambda_{M}} \\
M'
\ar[r]^(0.4){g^{\sharp}(\rho_M)} 
& \opn{Sq}_{B' / A} M'
} \]
of morphisms in $\msf{D}(\cat{Mod} B)$ is commutative.
\end{lem}

\begin{proof}
Recall that $g^{\sharp}(\rho_M)$ was defined to be
$g^{\sharp}(\rho^{\mrm{tau}}_B) \otimes \rho_M$,
under the isomorphism
\[ (\opn{Sq}_{B / A} M) \otimes^{\mrm{L}}_{B} (\opn{Sq}_{B' / B} B')
\xar{\smallsmile_{g^*; M, B'}}
\opn{Sq}_{B' / A} M' . \]
But in this situation of a localization 
$B' \otimes_B B' = B'$, and 
\[ g^{\sharp}(\rho^{\mrm{tau}}_B) : B' \iso
\opn{Sq}_{B' / B} B' =
\opn{RHom}_{B' \otimes_B B'}(B', B' \otimes_B B')  = B' \]
is just the identity. The composed map
\[ \opn{Sq}_{B / A} M \xar{\opn{q}}
(\opn{Sq}_{B / A} M) \otimes^{\mrm{L}}_{B} B'
\xar{\smallsmile_{g^*; M, B'}}
\opn{Sq}_{B' / A} M' \]
is exactly $\lambda_M$. 
\end{proof}

\begin{proof}[Proof of Theorem \tup{\ref{thm6.3}}]
Consider the cubic diagram
\[ \UseTips  \xymatrix @=8ex {
N 
\ar[r]^{\phi} 
\ar[d]_{\rho_N}
\ar[drrr]^(0.7){\mrm{q}}
& M 
\ar[d]|(0.39)\hole^(0.6){\rho_M}
\ar[drrr]*++++{\hole}^{\mrm{q}}
\\
\opn{Sq}_{C / A} N
\ar[r]^(0.48){\opn{Sq}(\phi)} 
%{\opn{Sq}_{f^* / A}(\phi) } 
\ar[drrr]^(0.5){\lambda_N}
& \opn{Sq}_{B / A} M
\ar[drrr]*+++++++{\hole}|(0.59)\hole^(0.3){\lambda_M}
& & N'
\ar[r]_{\phi'}
\ar[d]^(0.4){h^{\sharp}(\rho_N)}
& M' 
\ar[d]_(0.4){g^{\sharp}(\rho_M)} \\
& & & \opn{Sq}_{C' / A} N'
\ar[r]_(0.47){\opn{Sq}(\phi')} 
%{\opn{Sq}_{{f'}^* / A}(\phi')}
& \opn{Sq}_{B' / A} M'
} \]
The top face of the diagram is trivially commutative. 
The right face is commutative by Lemma \ref{lem6.5}. 
For the same reason the left face 
commutes. The bottom face commutes by Lemma \ref{lem6.4}.
Because $\rho_M$ and $\rho_N$ are isomorphisms we see that
$B' \otimes_B \opn{Sq}_{C / A} N \cong \opn{Sq}_{C' / A} N'$
and
$B' \otimes_B \opn{Sq}_{B / A} M \cong \opn{Sq}_{B' / A} M'$,
via the morphisms $\lambda_N$ and $\lambda_M$ respectively.
Thus the front face is gotten from the rear face by applying the 
functor $g^{\sharp} = B' \otimes_B -$ to it. 
But the rear face is commutative because $\phi$ is a rigid trace 
morphism. We conclude that the front face is commutative, which 
means that $\phi'$ is a rigid trace morphism. 
\end{proof}

\begin{rem} \label{rem6.2}
Theorem \ref{thm6.3} is true even when $g^*$ is essentially 
\'etale (not just a localization). 
However the proof is more complicated, since in this 
case $B' \otimes_B B'$ is not equal to $B'$; cf.\ Proposition 
\ref{prop3.5}. In the definition of the morphism $\lambda_M$, 
namely in equation (\ref{eqn6.3}), one has to replace the 
isomorphism $\diamondsuit$ with 
\[ \begin{aligned}
\mrm{RHom}_{\til{B} \otimes_{A} \til{B}}(B, 
M' \otimes^{\mrm{L}}_{A} M') & \cong
\mrm{RHom}_{\til{B}' \otimes_{A} \til{B'}}(B' \otimes_B B',
M' \otimes^{\mrm{L}}_{A} M') \\
& \xar{(\nu, \bsym{1} \otimes \bsym{1})}
\mrm{RHom}_{\til{B}' \otimes_{A} \til{B'}}(B',
M' \otimes^{\mrm{L}}_{A} M') ,
\end{aligned} \]
where $\nu : B' \to B' \otimes_B B'$ is the canonical homomorphism 
from Proposition \ref{prop3.5}.
\end{rem}

We end the paper with an existence result. 

\begin{thm} \label{thm6.4}
Let $\K$ be a regular noetherian ring of finite Krull 
dimension, and let $A$ be an essentially finite type $\K$-algebra.
Then $A$ has a rigid complex $(R, \rho)$ relative to $\K$, 
such that the canonical homomorphism
$A \to \opn{Hom}_{\msf{D}(\cat{Mod} A)}(R, R)$
is bijective. 
\end{thm}

\begin{proof}
First we note that 
$\opn{Sq}_{\K / \K} \K = \K$, and this gives the tautological 
rigid complex 
$(\K, \rho^{\mrm{tau}}) \in
\msf{D}^{\mrm{b}}_{\mrm{f}}(\cat{Mod} \K)_{\mrm{rig} / \K}$.

Now the structural homomorphism $\K \to A$ can be factored into
\[ \K \xar{f^*} B \xar{g^*} C \xar{h^*} A , \]
where $B := \K[t_1, \ldots, t_n]$ is a polynomial algebra;
$B \to C$ is surjective; and $C \to A$ is a localization. 
By definition 
$f^{\sharp} \K = \Omega^n_{B / \K}[n]
\in \msf{D}^{\mrm{b}}_{\mrm{f}}(\cat{Mod} B)$,
and according to Theorem \ref{thm2.4} it has an induced 
rigidifying isomorphism $f^{\sharp}(\rho^{\mrm{tau}})$.
Since $B$ is a regular ring of finite Krull dimension it follows 
that 
\[ g^{\flat} f^{\sharp} \K = \opn{RHom}_B(C, 
\Omega^n_{B / \K}[n]) 
\in \msf{D}^{\mrm{b}}_{\mrm{f}}(\cat{Mod} C) . \]
Since $\K$ is regular it follows that 
$g^{\flat} f^{\sharp} \K$ has finite flat dimension over $\K$.
By Theorem \ref{thm2.2} we deduce that 
$g^{\flat} f^{\sharp} \K$ has an induced rigidifying isomorphism 
\linebreak
$g^{\flat} f^{\sharp} (\rho^{\mrm{tau}})$.
Using Theorem \ref{thm2.4} again we obtain a rigidifying 
isomorphism
$\rho := h^{\sharp} g^{\flat} f^{\sharp} (\rho^{\mrm{tau}})$
for the complex
$R := h^{\sharp} g^{\flat} f^{\sharp} \K
\in \msf{D}^{\mrm{b}}_{\mrm{f}}(\cat{Mod} A)$.

It remains to calculate 
$\opn{RHom}_A(R, R)$. Since $\Omega^n_{B / \K}$ is a projective 
$B$-module of rank $1$ we see that 
\[ \opn{RHom}_B(f^{\sharp} \K, f^{\sharp} \K) \cong B . \]
Next, by twisting, adjunction and the fact that $B$ is regular
we get
\[ \begin{aligned}
& \opn{RHom}_C(g^{\flat} f^{\sharp} \K, 
g^{\flat} f^{\sharp} \K) 
\cong \opn{RHom}_C(g^{\flat} B, g^{\flat}B) \\
& \quad =
\opn{RHom}_C \bigl( \opn{RHom}_B(C, B), \opn{RHom}_B(C, B) 
\bigr) \\
& \quad \cong \opn{RHom}_B \bigl( \opn{RHom}_B(C, B), B \bigr) 
\cong C . 
\end{aligned} \]
Finally, using Proposition \ref{prop2.4}, and adjunction, we get
\[ \begin{aligned}
& \opn{RHom}_A(h^{\sharp}  g^{\flat} f^{\sharp} \K, 
h^{\sharp}  g^{\flat} f^{\sharp} \K) 
\cong \opn{RHom}_C(g^{\flat} f^{\sharp} \K, 
h^{\sharp}  g^{\flat} f^{\sharp} \K) \\
& \quad \cong
A \otimes_C \opn{RHom}_C(g^{\flat} f^{\sharp} \K, 
g^{\flat} f^{\sharp} \K) 
\cong A \otimes_C C \cong A . 
\end{aligned} \]
\end{proof}

%\newpage

\end{document}